\documentclass[11pt,leqno]{smfart}

\usepackage[latin1]{inputenc}
\usepackage[frenchb]{babel}
\usepackage{amsmath}
\usepackage{amssymb}
\usepackage{a4wide}
\usepackage{epsfig}
\usepackage{stmaryrd}

\newtheorem{theo}{Th\'eor\`eme}[section]
\newtheorem{lemm}[theo]{Lemme}
\newtheorem{prop}[theo]{Proposition}

\newtheorem{rema}[theo]{Remarque}

\newcounter{hypo}
\newcounter{hyphyp}

\newenvironment{hyp}{\renewcommand{\theenumi}{{\rm (H\arabic{enumi})}} \renewcommand{\labelenumi}{{\bf (H\arabic{enumi})}} \begin{enumerate} \setcounter{enumi}{\value{hypo}} \item}{\stepcounter{hypo} \end{enumerate}}

\newenvironment{hypri}{\addtocounter{hypo}{-1} \setcounter{hyphyp}{\value{hypo}} \renewcommand{\theenumi}{{\rm (H\arabic{enumi})'}} \renewcommand{\labelenumi}{{\bf (H\arabic{enumi})'}} \begin{enumerate} \setcounter{enumi}{\value{hyphyp}} \item}{\stepcounter{hypo} \end{enumerate}}

\def\C{{\mathbb C}}
\def\N{{\mathbb N}} 
\def\R{{\mathbb R}}

\def\S{{\mathbb S}}

\def\CI{\mathcal {I}}

\def\CO{\mathcal {O}}

\def\re{\mathop{\rm Re}\nolimits}
 \def\im{\mathop{\rm Im}\nolimits}

\def\Op{\mathop{\rm Op}\nolimits}

\def\supp{\mathop {\rm supp}\nolimits}

\def\sgn{\mathop{\rm sgn}\nolimits}
\def\Hess{\mathop{\rm Hess}\nolimits}
\def\<{\langle}
\def\>{\rangle}

\def\dl{\mathcal{D}}

\newcommand{\fract}[2]{\genfrac{}{}{0pt}{}{\scriptstyle #1}{\scriptstyle #2}}

\renewcommand{\theenumi}{\sl{\roman{enumi}}}
\renewcommand{\labelenumi}{\theenumi)}

\makeatletter
 \@addtoreset{equation}{section}
 \makeatother
 \renewcommand{\theequation}{\thesection.\arabic{equation}}

\author{Jean-Fran\c{c}ois Bony}

\address{\newline Institut de Math\'ematiques de Bordeaux, \newline UMR 5251 du CNRS, \newline Universit\'e de Bordeaux I, \newline 351 cours de la Lib\'eration, \newline 33405 Talence, \newline France \newline}

\email{bony@math.u-bordeaux1.fr}

\title{Mesures limites pour l'\'equation de Helmholtz dans le cas non captif}

\subjclass{35Q40, 35J10, 35S30, 81Q20}

\keywords{Equation de Helmholtz, Mesure semi-classique, Analyse micro-locale}

\begin{document}

\begin{abstract}
Cet article est consacr\'e \`a l'\'etude des mesures limites associ\'ees \`a la solution de l'\'equation de Helmholtz avec un terme source se concentrant en un point. Le potentiel est suppos\'e $C^{\infty}$ et l'op\'erateur non-captif. La solution de l'\'equation de Schr\"odinger semi-classique s'\'ecrit alors micro-localement comme somme finie de distributions lagrangiennes.

Sous une hypoth\`ese g\'eom\'etrique, qui g\'en\'eralise l'hypoth\`ese du viriel, on en d\'eduit que la mesure limite existe et qu'elle v\'erifie des propri\'et\'es standard. Enfin, on donne un exemple d'op\'erateur qui ne v\'erifie pas l'hypoth\`ese g\'eom\'etrique et pour lequel la mesure limite n'est pas unique. Le cas de deux termes sources est aussi trait\'e.
\end{abstract}

\maketitle


\section{Introduction.}

Dans ce papier, on \'etudie la limite hautes fr\'equences de l'\'equation de Helmholtz, en dimension $n \geq 1$, avec un terme source se concentrant en $x=0$. Par un changement de variables, ce probl\`eme se ram\`ene \`a l'\'etude de la limite quand $h  \to 0$ de $u_{h}$, la solution de
\begin{equation*}
\Big( - \frac{h^{2}}{2} \Delta + V(x) -E \Big) u_{h} = h^{- n/2} S \Big( \frac{x}{h} \Big) .
\end{equation*}
Le membre de droite de cette \'equation est une fonction qui joue le r\^ole d'un profil se concentrant en $0$ \`a l'\'echelle $h$. Sa norme $L^{2}$ est pr\'eserv\'ee. L'\'energie, ou param\`etre de r\'egularisation, $E$ est dans le demi-plan complexe sup\'erieur et v\'erifie $E = E_{0} + h E_{1} + o (h)$, $E_{0}>0$, quand $h \to 0$. De plus le comportement du potentiel $V$ est prescrit \`a l'infini (dans cet article, $V$ sera \`a longue port\'ee).

Plus pr\'ecis\'ement, on veut calculer les \'eventuelles mesures semi-classiques de Wigner associ\'ees \`a $u_{h}$. Autrement dit, on cherche \`a savoir si, pour $q \in C^{\infty}_{0} ( \R^{2n})$, la quantit\'e
\begin{equation*}
\lim_{h \to 0} \, h \< \Op (q) u_{h} , u_{h} \> ,
\end{equation*}
existe et s'\'ecrit $\< q, \mu_{h^{-1}} \>$ o\`u $\mu_{h^{-1}}$ est une mesure sur ${\rm T}^{\star} \R^{n} = \R^{2n}$ ind\'ependante de $h$ ($\Op (q)$ est l'op\'erateur pseudo-diff\'erentiel d\'efini en \eqref{a78}). La pr\'esence du $h$ devant le produit scalaire sera obligatoire pour rendre la quantit\'e convergente.

En g\'en\'eral cette mesure v\'erifie les propri\'et\'es suivantes. Elle est support\'ee dans la surface d'\'energie $E_{0}$. Elle est nulle dans la zone entrante (on parle aussi de condition de radiation \`a l'infini). Enfin cette mesure limite v\'erifie l'\'equation de Liouville :
\begin{equation*}
 \big( \xi \partial_{x} - \nabla V (x) \partial_{\xi} + 2 \im E_{1} \big) \mu_{h^{-1}} = (2 \pi)^{1-n} \vert \widehat{S} (\xi ) \vert^{2} \delta_{x =0} \delta_{\frac{1}{2} \xi^{2} + V(0) = E_{0}} .
\end{equation*}
On verra que le dernier point n'est pas toujours satisfait.

A notre connaissance, le premier travail dans cette direction est du \`a Benamou, Castella, Katsaounis et Perthame \cite{BeCaKaPe02_01}. Sous l'hypoth\`ese que $V =0$, ils obtiennent la convergence vers la mesure pr\'ec\'edemment d\'ecrite. Ce r\'esultat a \'et\'e \'etendu au cas des potentiels variables par Wang \cite{Wa07_01}. Lorsque la source est port\'ee par une vari\'et\'e, des r\'esultats \'equivalents ont \'et\'e obtenus par Castella, Perthame et Runborg \cite{CaPeRu02_01} dans le cas $V =0$ et g\'en\'eralis\'es par Wang et Zhang \cite{WaZh06_01} aux cas \`a coefficients variables. Enfin, quand le potentiel est constant de chaque cot\'e d'un hyperplan sur lequel il pr\'esente une discontinuit\'e, la propagation de telles mesures semi-classiques a \'et\'e trait\'ee par Fouassier \cite{Fo07_01}. Tous ces travaux utilisent des majorations de la r\'esolvante entre espaces de Besov (estimations de Morrey--Campanato). Ces estimations sont parfois obtenues \`a l'aide de la th\'eorie de Mourre. C'est pourquoi la r\'egularit\'e du potentiel peut \^etre abaiss\'ee jusqu'\`a $C^{2}$, dans certains cas. N\'eanmoins, tous supposent que le potentiel est nul ou satisfait l'hypoth\`ese du viriel. Ce pr\'esent travail a pour but de d\'ecrire ce qui se passe quand cette hypoth\`ese est retir\'ee, mais que l'op\'erateur reste non captif.

Mentionnons qu'il existe \'egalement des articles qui prouvent sp\'ecifiquement des estimations de la r\'esolvante entre espaces de Besov (voir, par exemple, Perthame et Vega \cite{PeVe99_01} pour des potentiels ne s'annulant pas \`a l'infini, Wang et Zhang \cite{WaZh06_01} pour des termes sources port\'es par des surfaces, Castella et Jecko \cite{CaJe06_01} dans le cas $C^{2}$ non-captif, Wang \cite{Wa07_01} et Castella, Jecko et Knauf \cite{CaJeKn07_01} pour des singularit\'es coulombiennes). Enfin, Castella \cite{Ca05_01} a \'etudi\'e la limite faible hautes fr\'equences de l'\'equation de Helmholtz dans le cas non captif.

Dans cet article, on d\'emontre que la solution de l'\'equation de Schr\"odinger stationnaire peut s'\'ecrire micro-localement pr\`es de tout point $\rho_{0} \in \R^{2 n}$ comme somme finie de distributions lagrangiennes. Toutefois, la vari\'et\'e $x=0$ pose des difficult\'es, et on peut juste majorer la solution pr\`es de l'origine. Une fois cette construction effectu\'ee, il est facile de d\'eduire l'existence de la mesure limite et de la calculer. A cause du caract\`ere bilin\'eaire de ces mesures, on doit faire une hypoth\`ese g\'eom\'etrique sur les vari\'et\'es lagrangiennes qui portent la solution. Cette hypoth\`ese est plus g\'en\'erale que l'hypoth\`ese du viriel. Elle, ou une autre du m\^eme type, est en fait obligatoire : dans la partie \ref{a50}, en s'inspirant de l'exemple de Castella \cite[section 9]{Ca05_01}, on construit un op\'erateur qui ne v\'erifie pas l'hypoth\`ese g\'eom\'etrique et pour lequel la mesure limite n'est pas unique.

Pour obtenir ces r\'esultats, on utilise une m\'ethode d\'ependante du temps et on \'ecrit la r\'esolvante comme l'int\'egrale du propagateur en temps positifs. Des constructions BKW d\'evelopp\'ees par Maslov et H\"ormander servent alors \`a approcher le propagateur. Cette id\'ee avait \'et\'e utilis\'ee par Castella dans \cite{Ca05_01}, mais elle vient surtout, pour ce type de questions, de l'\'etude de l'amplitude de diffusion r\'ealis\'ee par Robert et Tamura \cite{RoTa89_01} en 1989. Remarquons que, en un certain sens, l'amplitude de diffusion peut \^etre vue comme l'\'equation de Helmholtz avec terme source \`a l'infini. Ces m\'ethodes, tr\`es gourmandes en d\'eriv\'ees, n\'ecessitent que le potentiel et la transform\'ee de Fourier du terme source soient essentiellement $C^{\infty}$. Par contre, la preuve n'utilise ni la th\'eorie de Mourre ni les estimations de la r\'esolvante entre espaces de Besov. Cependant le lemme \ref{a7}, qui sert \`a contr\^oler la norme $L^{2}$ (d'une partie) de la solution tr\`es pr\`es de $x=0$ utilise une d\'ecomposition dyadique.

Enfin, avec cette m\'ethode, on peut aussi traiter le cas de deux (ou plusieurs) points sources. Nous aurons besoin d'une hypoth\`ese g\'eom\'etrique qui, comme pr\'ec\'edemment, est essentiellement n\'ecessaire pour l'unicit\'e de la mesure. Fouassier \cite{Fo06_01} avait \'etudi\'e cette question sous l'hypoth\`ese $V =0$.

On se place dans le cadre de l'analyse semi-classique et $h \in ]0,h_{0}[$ d\'esigne le petit param\`etre. Dans la suite, on utilisera les notations suivantes. Pour $m (x, \xi) \geq 0$ une fonction de poids, $S (m)$ d\'esigne l'ensemble des fonctions $q (x, \xi ,h) \in C^{\infty} (\R^{2n}_{x, \xi})$ telles que
\begin{equation}
\partial_{x}^{\alpha} \partial_{\xi}^{\beta} q (x, \xi ,h) = \CO (m),
\end{equation}
pour tous $\alpha , \beta \in \N^{n}$. On dit que le symbole $q (x, \xi ,h)$ est classique et on note $q \in S_{{\rm cl}} (m)$ si il existe une suite de symboles $q_{j} ( x, \xi ) \in S (m)$ tel que
\begin{equation}
q (x, \xi ,h) - \sum_{j=0}^{N} q_{j} (x, \xi ) h^{j} \in S (h^{N+1} m) :=h^{N+1} S (m) .
\end{equation}
Sous cette hypoth\`ese, $q_{0}$ est appel\'e le symbole principal. Pour $q \in S (m)$, l'op\'erateur pseudo-diff\'erentiel (semi-classique en quantification de Weyl) de symbole $q$ est d\'efini par
\begin{equation} \label{a78}
\Op (q) f (x) = \frac{1}{(2 \pi h)^{n}} \iint e^{i (x - y) \cdot \xi /h} q \Big( \frac{x- y}{2} , \xi ,h \Big) f (y) \, d y \, d \xi ,
\end{equation}
au sens des int\'egrales oscillantes. On note $\Psi (m) = \Op ( S (m))$, l'ensemble des op\'erateurs pseudo-diff\'erentiels de symbole dans $S (m)$. Enfin,
\begin{equation}
\widehat{f} (\xi) = \frac{1}{( 2 \pi )^{n/2}} \int e^{- i x \cdot \xi /h} f(x) \, d x ,
\end{equation}
est la transform\'ee de Fourier de $f \in {\mathcal S}' (\R^{n})$, le dual de l'espace de Schwartz ${\mathcal S} (\R^{n})$.

\section{Enonc\'e des r\'esultats.}

On consid\`ere un op\'erateur de Schr\"odinger semi-classique sur $L^{2}(\R^{n})$, avec $n \geq 1$,
\begin{equation}
P = - \frac{h^{2}}{2} \Delta + V(x) ,
\end{equation}
et on fait les hypoth\`eses suivantes.

\begin{hyp} \label{h1}
$V$ est une fonction $C^\infty$ sur $\R^{n}$, et il existe $\rho> 0$ tel que
\begin{equation*}
\vert \partial^{\alpha} V (x) \vert \lesssim \< x \>^{- \rho - \vert \alpha \vert} ,
\end{equation*}
pour tout $\alpha \in \N^{n}$.
\end{hyp}
Autrement dit le potentiel est r\'egulier et \`a longue port\'ee. Le symbole de $P$ est not\'e
\begin{equation*}
p (x, \xi ) = \frac{\xi^{2}}{2} + V(x) ,
\end{equation*}
et son champ hamiltonien
\begin{equation*}
{\rm H}_{p} = \xi \frac{\partial \ }{\partial x} - \nabla V (x) \frac{\partial \ }{\partial \xi} .
\end{equation*}

\begin{hyp} \label{h2}
La surface d'\'energie $E_{0} >0$ est non captive. C'est \`a dire,
\begin{equation*}
{\mathcal K}(E_{0} ) = \{ (x, \xi ) \in p^{-1} (E_{0} ) ; \ t \to \exp (t {\rm H}_{p}) ( x, \xi ) \text{ reste born\'e en temps}\} = \emptyset.
\end{equation*}
\end{hyp}

\begin{hyp} \label{h3}
De plus, $V (0) < E_{0}$.
\end{hyp}

\begin{figure}
\begin{center}
\begin{picture}(0,0)%
\includegraphics{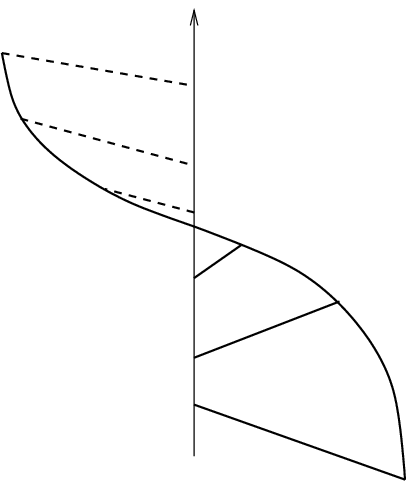}%
\end{picture}%
\setlength{\unitlength}{1184sp}%
\begingroup\makeatletter\ifx\SetFigFont\undefined%
\gdef\SetFigFont#1#2#3#4#5{%
  \reset@font\fontsize{#1}{#2pt}%
  \fontfamily{#3}\fontseries{#4}\fontshape{#5}%
  \selectfont}%
\fi\endgroup%
\begin{picture}(6516,7813)(2893,-6952)
\put(7426,-1036){\makebox(0,0)[lb]{\smash{{\SetFigFont{9}{10.8}{\rmdefault}{\mddefault}{\updefault}$\R^{2 n}$}}}}
\put(7801,-5236){\makebox(0,0)[lb]{\smash{{\SetFigFont{9}{10.8}{\rmdefault}{\mddefault}{\updefault}$\widetilde{\Lambda}$}}}}
\put(5626,-6811){\makebox(0,0)[lb]{\smash{{\SetFigFont{9}{10.8}{\rmdefault}{\mddefault}{\updefault}$x=0$}}}}
\end{picture} \qquad \qquad \qquad
\begin{picture}(0,0)%
\includegraphics{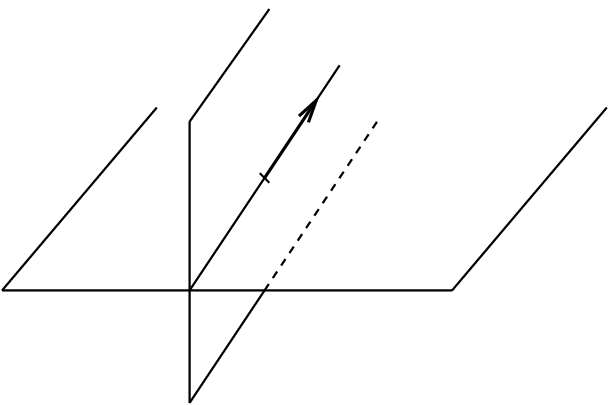}%
\end{picture}%
\setlength{\unitlength}{1184sp}%
\begingroup\makeatletter\ifx\SetFigFont\undefined%
\gdef\SetFigFont#1#2#3#4#5{%
  \reset@font\fontsize{#1}{#2pt}%
  \fontfamily{#3}\fontseries{#4}\fontshape{#5}%
  \selectfont}%
\fi\endgroup%
\begin{picture}(9741,6366)(1768,-6394)
\put(7501,-5836){\makebox(0,0)[lb]{\smash{{\SetFigFont{9}{10.8}{\rmdefault}{\mddefault}{\updefault}$\R^{2 n}$}}}}
\put(8326,-4261){\makebox(0,0)[lb]{\smash{{\SetFigFont{9}{10.8}{\rmdefault}{\mddefault}{\updefault}$\Lambda_{1}$}}}}
\put(5101,-2011){\makebox(0,0)[lb]{\smash{{\SetFigFont{9}{10.8}{\rmdefault}{\mddefault}{\updefault}$\Lambda_{2}$}}}}
\put(6151,-3061){\makebox(0,0)[lb]{\smash{{\SetFigFont{9}{10.8}{\rmdefault}{\mddefault}{\updefault}$\rho_{0}$}}}}
\put(6901,-2011){\makebox(0,0)[lb]{\smash{{\SetFigFont{9}{10.8}{\rmdefault}{\mddefault}{\updefault}${\rm H}_{p}$}}}}
\put(4126,-5161){\makebox(0,0)[lb]{\smash{{\SetFigFont{9}{10.8}{\rmdefault}{\mddefault}{\updefault}$\Lambda$}}}}
\end{picture}%
\caption{La vari\'et\'e $\widetilde{\Lambda}$ et la situation g\'eom\'etrique pr\`es d'un point $\rho_{0}$ o\`u $\Lambda$ se ``recoupe''.} \label{a66}
\end{center}
\end{figure}

\begin{hyp}  \label{h5}
On suppose que
\begin{equation*}
\text{mes}_{n-1} \big\{ \xi \in \sqrt{2 (E_{0} - V(0))} \S^{n-1} ; \ \exists t > 0 \quad \Pi_{x} \exp ( t {\rm H}_{p}) (0 , \xi) = 0 \big\} =0 .
\end{equation*}
$\text{mes}_{n-1}$ est la mesure de Lebesgue sur $\sqrt{2 (E_{0} - V(0))} \S^{n-1}$ et $\Pi_{x} ( x, \xi ) = x$ est la projection spatiale.
\end{hyp}

\begin{figure}
\begin{center}
\begin{picture}(0,0)%
\includegraphics{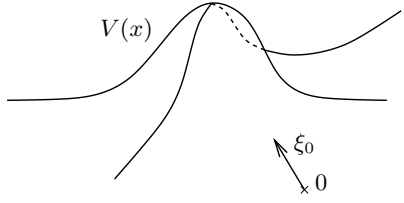}%
\end{picture}%
\setlength{\unitlength}{1184sp}%
\begingroup\makeatletter\ifx\SetFigFont\undefined%
\gdef\SetFigFont#1#2#3#4#5{%
  \reset@font\fontsize{#1}{#2pt}%
  \fontfamily{#3}\fontseries{#4}\fontshape{#5}%
  \selectfont}%
\fi\endgroup%
\begin{picture}(8316,4087)(1768,-4777)
\put(3751,-1411){\makebox(0,0)[lb]{\smash{{\SetFigFont{9}{10.8}{\rmdefault}{\mddefault}{\updefault}$V(x)$}}}}
\put(8251,-4636){\makebox(0,0)[lb]{\smash{{\SetFigFont{9}{10.8}{\rmdefault}{\mddefault}{\updefault}$0$}}}}
\put(7801,-3811){\makebox(0,0)[lb]{\smash{{\SetFigFont{9}{10.8}{\rmdefault}{\mddefault}{\updefault}$\xi_{0}$}}}}
\end{picture}%
\caption{Un exemple de potentiel $V$ qui, pour $0 < E_{0} < \max \, V$, ne v\'erifie pas l'hypoth\`ese du viriel mais pour lequel \ref{h5} est vrai. En effet, l'ensemble qui apparait dans \ref{h5} est un singleton $\{ \xi_{0} \}$.} \label{a83}
\end{center}
\end{figure}

\noindent
Notons que l'ensemble qui apparait dans \ref{h5} est mesurable car ferm\'e. L'interpr\'etation g\'eom\'etrique de cette hypoth\`ese est la suivante (voir la partie \ref{a85} pour les preuves). La solution de l'\'equation de Helmholtz sera micro-localis\'ee sur l'ensemble $\Lambda$ d\'efini par
\begin{equation}
\Lambda := \{ \exp ( t {\rm H}_{p} ) ( 0, \xi) ; \ t >0 \text{ et } \xi^{2} = 2 (E_{0} - V (0)) \} \subset p^{-1}( E_{0} ) .
\end{equation}
Cet ensemble peut se ``recouper'' (voir les figures \ref{a66}, \ref{a83}, \ref{a57} (en dimension $n=1$) et \ref{a56} (en dimension sup\'erieure)). Mais localement pr\`es d'un point $\rho_{0} = ( x_{0} , \xi_{0} ) \in p^{-1} (E_{0} )$, il est de la forme
\begin{equation}
\Lambda = \bigcup_{k=1}^{K} \Lambda_{k} , \text{ si } x_{0} \neq 0 \quad \text{ et } \quad \Lambda = \widetilde{\Lambda} \cup \bigcup_{k=1}^{K} \Lambda_{k}, \text{ si } x_{0} =0 ,
\end{equation}
o\`u les ensembles $\widetilde{\Lambda}$ et $\Lambda_{k}$ sont d\'efinis ci-dessous. La vari\'et\'e
\begin{equation} \label{a79}
\widetilde{\Lambda} = \{ \exp ( t {\rm H}_{p} ) ( 0, \xi) ; \ 0<t <\varepsilon_{0} \text{ et } \xi^{2} = 2 (E_{0} - V (0)) \} ,
\end{equation}
avec $\varepsilon_{0} >0$ assez petit, est une vari\'et\'e (locale) lagrangienne $C^{\infty}$ se projetant bien suivant les variables d'espace, pour $x \neq 0$. Elle est repr\'esent\'ee sur la figure \ref{a66}. Les $\Lambda_{k}$ sont des vari\'et\'es (locales) lagrangiennes $C^{\infty}$. Plus pr\'ecis\'ement,
\begin{equation}  \label{a69}
\{ (t, \xi) ; \ 0<t , \ \xi^{2} = 2 (E_{0} - V (0))  \text{ et }\exp ( t {\rm H}_{p}) (0 , \xi) = \rho_{0} \} ,
\end{equation}
est un ensemble fini $\{ (t_{1} , \xi_{1} ) , \ldots , (t_{K} , \xi_{K} ) \}$ et
\begin{equation}  \label{a53}
\Lambda_{k} = \{ \exp ( t {\rm H}_{p} ) ( 0, \xi) ; \ t \text{ au voisinage de } t_{k} \text{ et } \xi \text{ au voisinage de } \xi_{k} \} .
\end{equation}
Ainsi, il se peut que $\Lambda_{j} = \Lambda_{k}$ avec $j \neq k$. On notera parfois $\Lambda_{0} = \widetilde{\Lambda}$ et, d'apr\`es la section \ref{a85}, l'hypoth\`ese \ref{h5} est \'equivalente \`a

\begin{hypri} \label{h7}
Les intersections des vari\'et\'es $\widetilde{\Lambda}$ et $\Lambda_{k}$ sont de mesure nulle :
\begin{equation*}
\text{mes}_{n} ( \Lambda_{j} \cap \Lambda_{k} ) =0 ,
\end{equation*}
pour $0 \leq j \neq k \leq K$. La mesure est prise sur la vari\'et\'e $\Lambda_{j}$ (ou $\Lambda_{k}$).
\end{hypri}

\noindent
D'apr\`es le principe d'absorption limite, on sait que, pour tout $E \in ] 0 , + \infty [$,
\begin{equation}
(P- (E+ i 0))^{-1} := \lim_{\delta \to 0} (P - (E + i \delta))^{-1} ,
\end{equation}
existe en tant qu'op\'erateur de $L^{2} ( \< x \>^{2 \alpha} d x)$ dans $L^{2} ( \< x \>^{- 2 \alpha} d x)$, pour tout $\alpha > 1/2$. De plus, Robert et Tamura \cite{RoTa87_01} ont d\'emontr\'e que, sous les hypoth\`eses \ref{h1} et \ref{h2},
\begin{equation}  \label{a1}
\sup_{0 \leq \delta < 1} \big\Vert \<x \>^{- \alpha} (P - (E + i \delta))^{-1} \< x \>^{- \alpha} \big\Vert \lesssim h^{-1} ,
\end{equation}
pour $E \in ]0 , + \infty [$ au voisinage de $E_{0}$. Leur preuve utilise, entre autres choses, les constructions d'Isozaki et Kitada. Mais il existe d'autres d\'emonstrations, notamment celle de C. G\'erard et Martinez \cite{GeMa88_01} qui est bas\'ee sur la m\'ethode de Mourre \cite{Mo81_01}. D'autre part, Burq \cite{Bu02_01} a prouv\'e \eqref{a1} par l'absurde grâce aux mesures de d\'efauts. Cette estimation r\'esulte aussi du r\'ecent papier de C. G\'erard \cite{Ge07_01} qui en donne une preuve directe utilisant uniquement de l'analyse fonctionnelle. Dans la suite et dans ce genre d'estimation, $\alpha$ d\'esignera toujours un r\'eel strictement sup\'erieur \`a $1/2$.

On prend une donn\'ee $S_{h}$ qui se concentre en $0$ de la mani\`ere suivante :
\begin{hyp} \label{h4}
On consid\`ere la suite de fonctions
\begin{equation*}
S_{h} (x) = h^{-n/2} S \Big( \frac{x}{h} \Big) ,
\end{equation*}
o\`u $S$ est dans $\< x \>^{-N} L^{2} (\R^{n})$ pour tout $N \geq 0$. Donc, $\widehat{S} \in C^{\infty} (\R^{n} )$.
\end{hyp}
et on fait l'hypoth\`ese suivante sur l'\'energie :
\begin{hyp} \label{h6}
Soit $E = E (h) \in \C$ une famille d'\'energies v\'erifiant
\begin{equation*}
E = E_{0} + h E_{1} + o (h) ,
\end{equation*}
et $\im E \geq 0$.
\end{hyp}
Ce papier est consacr\'e \`a l'\'etude de la solution de Helmholtz
\begin{equation}
u_{h} (x) = (P- (E+ i 0))^{-1} S_{h} .
\end{equation}
Si $\im E >0$, on a $u_{h} = (P - E)^{-1} S_{h}$.

\begin{theo}\sl \label{a48}
Sous les hypoth\`eses \ref{h1}--\ref{h6}, il existe une mesure de Radon positive $\mu_{h^{-1}}$ sur $\R^{2n}$ telle que, pour tout $q \in C^{\infty}_{0} (\R^{2n})$,
\begin{equation}  \label{a68}
h \< \Op (q) u_{h} , u_{h} \> \longrightarrow \< q , \mu_{h^{-1}} \> ,
\end{equation}
quand $h \to 0$. La mesure $\mu_{h^{-1}}$ ne d\'epend pas de $h$ et l'indice $h^{-1}$ est juste une notation.

\noindent
La mesure $\mu_{h^{-1}}$ poss\`ede les propri\'et\'es suivantes qui la caract\'erisent uniquement :
\begin{enumerate}
\item La mesure $\mu_{h^{-1}}$ est support\'ee dans la surface d'\'energie $E_{0}$ : \label{a80}
\begin{equation}
\supp \mu_{h^{-1}} \subset p^{-1}(E_{0}) .
\end{equation}

\item La mesure $\mu_{h^{-1}}$ est nulle dans la zone entrante : \label{a81}
\begin{equation}
\mu_{h^{-1}} =0 \text{ dans } \{ (x, \xi ) ; \ \vert x \vert \text{ grand  et } \cos (x, \xi ) < -1/2 \}.
\end{equation}

\item La mesure $\mu_{h^{-1}}$ v\'erifie l'\'equation d'\'evolution :\label{a82}
\begin{align}
( {\rm H}_{p} + 2 \im E_{1} ) \, \mu_{h^{-1}} =&  \big( \xi \partial_{x} - \nabla V (x) \partial_{\xi} + 2 \im E_{1} \big) \mu_{h^{-1}}   \nonumber  \\
=& (2 \pi)^{1-n} \vert \widehat{S} (\xi ) \vert^{2} \delta_{x =0} \delta_{\frac{1}{2} \xi^{2} + V(0) = E_{0}} .
\end{align}
\end{enumerate}
\end{theo}

Le membre de gauche de l'\'equation \eqref{a68} est \`a prendre au sens de la dualit\'e $L^{2} ( \< x \>^{2 \alpha} d x)$, $L^{2} ( \< x \>^{- 2 \alpha} d x)$. On dira qu'une fonction est un $\CO^{\alpha} ( 1 )$ si sa norme $L^{2} ( \< x \>^{- 2 \alpha} d x)$ est un $\CO (1)$. La mesure $\mu_{h^{-1}}$ est ind\'ependante de $h$. L'indice $h^{-1}$ sert juste \`a indiquer que le produit scalaire dans le membre de gauche de \eqref{a68} a \'et\'e multipli\'e par $h$.

\begin{rema}\sl
En fait, on calcule micro-localement $u_{h}$ pr\`es de tout point $\rho_{0} = ( x_{0} , \xi_{0} ) \in \R^{2n}$ avec $x_{0} \neq 0$ (pour $x_{0} =0$, on a juste les majorations donn\'ees par le lemme \ref{a7} et la discussion apr\`es le lemme \ref{a4} (voir aussi la proposition \ref{a59})). Si $p ( \rho_{0}) = E_{0}$ et $x_{0} \neq 0$, alors $u_{h}$ est somme de $K$ distributions lagrangiennes d'ordre $h^{-\frac{1}{2}}$, \`a symbole classique, port\'ees par les vari\'et\'es lagrangiennes $\Lambda_{k}$ (voir \eqref{a24} et \eqref{a54}). D'autre part, $u_{h}$ est nul micro-localement pr\`es des points $\rho_{0} \notin p^{-1} (E_{0})$ avec $x_{0} \neq 0$. Remarquons que ces r\'esultats ne n\'ecessitent pas l'hypoth\`ese \ref{h5}.

Comme corolaire, ces constructions donnent un d\'eveloppement asymptotique (au sens $C^{\infty}$) explicite de $u_{h}$ si les lagrangiennes $\Lambda_{k}$ se projettent bien spatialement. Enfin, la preuve du lemme \ref{a7} doit permettre d'\'etendre la construction de la distribution lagrangienne port\'ee par $\widetilde{\Lambda}$ jusqu'\`a des voisinages de taille $h^{1 - \varepsilon}$ de $0$ faisant apparaitre une phase singuli\`ere en $0$ comme dans le cas de $- h^{2} \Delta /2$.
\end{rema}

L'hypoth\`ese \ref{h5}, ou une hypoth\`ese du m\^eme type, est obligatoire pour obtenir les conclusions du th\'eor\`eme \ref{a48}. L'exemple donn\'e dans la partie \ref{a50} montre que, si cette hypoth\`ese est retir\'ee, la mesure limite n'est plus unique.

Dans la partie consacr\'ee \`a l'\'etude des mesures semi-classiques de l'\'equation de Schr\"odinger stationnaire de son papier \cite{Wa07_01}, Wang d\'emontre le th\'eor\`eme \ref{a48} en supposant que le potentiel $V \in C^{2} (\R^{n})$ satisfait
\begin{equation}
\nu_{0} (E_{0} - V(x)) - x \cdot \nabla V (x) \geq c_{0} ,
\end{equation}
avec $c_{0} >0$ et $\nu_{0} \in ]0, 2]$ (\'equation (4.23) de \cite{Wa07_01}). Comme $p = \xi^{2} + V(x)$ dans \cite{Wa07_01}, si $(x(t), \xi (t)) \in p^{-1} (E_{0} )$ est une courbe hamiltonienne, il vient
\begin{align}
\partial_{t} \big( x(t) \cdot \xi (t) \big) =& 2 \xi (t)^{2} - x \cdot \nabla V(x)  \nonumber  \\
\geq& \nu_{0} \xi (t)^{2} - x \cdot \nabla V(x)  \nonumber  \\
=& \nu_{0} (E_{0} - V(x)) - x \cdot \nabla V (x) \geq c_{0} .
\end{align}
En particulier, $x \cdot \xi$ est strictement croissant en temps et l'ensemble qui apparait dans l'hypoth\`ese \ref{h5} est vide.

L'hypoth\`ese \ref{h5} est proche de l'hypoth\`ese de transversalit\'e \'enonc\'ee par Castella dans \cite{Ca05_01}. Cet article est consacr\'e \`a la limite faible de l'\'equation de Helmholtz ; autrement dit, \`a la limite de $\< u_{h} , h^{- n/ 2} \varphi (\frac{x}{h} ) \>$, pour toute fonction $\varphi \in C^{\infty}_{0} (\R^{n})$. L'auteur suppose que la dimension $n$ est sup\'erieure \`a $2$, que l'ensemble
\begin{equation}
{\mathcal M} := \{ ( t , \xi , \eta ) ; \ \xi^{2} = 2 ( E_{0} - V(0) ) \text{ et } \exp ( t {\rm H}_{p}) ( 0 , \xi ) = ( 0 , \eta ) \} ,
\end{equation}
est une sous-vari\'et\'e de $\R^{ 2 n +1}$ de codimension strictement sup\'erieure \`a $n+2$ et que l'espace tangent de ${\mathcal M}$ est donn\'e par la lin\'earisation des \'equations qui le d\'efinissent (hypoth\`ese (H) page 245 de \cite{Ca05_01}). Mais alors, l'intersection de deux vari\'et\'es $\Lambda_{j}$ et $\Lambda_{k}$, pour $0 \leq j \neq k \leq K$, est une sous-vari\'et\'e de $\Lambda_{j}$ (ou $\Lambda_{k}$) de codimension sup\'erieure \`a $1$ et \ref{h7} est v\'erifi\'ee. N\'eanmoins, le probl\`eme trait\'e par Castella est diff\'erent de celui abord\'e ici et les comparaisons ne sont pas \'evidentes.

Le th\'eor\`eme \ref{a48} d\'ecrit la contribution principale de la solution de l'\'equation de Helmholtz. Mais il reste une partie sous-principale qui n'est pas principalement port\'ee par la surface d'\'energie $E_{0}$. Cette contribution est d\'ecrite par la proposition suivante.

\begin{prop}\sl \label{a59}
Sous les hypoth\`eses \ref{h1}--\ref{h6}, il existe une mesure de Radon positive $\mu_{1}$ sur $\R^{2n}$ telle que, pour tout $q \in C^{\infty}_{0} (\R^{2n})$,
\begin{equation}
\< \Op ( q (p-E_{0})^{2} ) u_{h} , u_{h} \> \longrightarrow \< q , \mu_{1} \> ,
\end{equation}
quand $h \to 0$. De plus, on a
\begin{equation}
\< q , \mu_{1} \> = \frac{1}{(2 \pi )^{n}} \int q (0, \xi ) \vert \widehat{S} ( \xi) \vert^{2} \, d \xi .
\end{equation}
\end{prop}

Le fait que $u_{h}$ porte une contribution particuli\`ere en $x=0$ est sp\'ecifique \`a l'\'etude de la r\'esolvante. Dans le cas de la mesure spectrale,
\begin{equation}
{\mathcal P} (E ) = \frac{1}{2 i \pi} \big( (P - (E + i 0))^{-1} - (P - (E - i 0))^{-1} \big) ,
\end{equation}
$E \in ] 0 , + \infty [$ proche de $E_{0}$, ce terme n'existe pas et ${\mathcal P} (E)$ est un op\'erateur int\'egral de Fourier r\'egulier. Ce r\'esultat a \'et\'e d\'emontr\'e par C. G\'erard et Martinez dans \cite{GeMa89_01}.

De plus, il semble aussi possible de d\'efinir une mesure $\mu_{h^{- 1/2}}$ par
\begin{equation}
h^{1/2} \< \Op ( q (p-E) ) u_{h} , u_{h} \> \longrightarrow \< q , \mu_{h^{-1/2}} \> ,
\end{equation}
pour tout $q \in C^{\infty}_{0} (\R^{2n})$. Sous une hypoth\`ese g\'eom\'etrique, cette mesure semble port\'ee par $p^{-1}(E_{0} ) \cap \{ x =0 \}$. En fait, cette mesure est un cas particulier de limite faible de l'\'equation de Helmholtz ; elle code le retour en $0$ de la vari\'et\'e $\Lambda$ lorsque l'hypoth\`ese de transversalit\'e (H) de Castella \cite{Ca05_01} n'est pas v\'erifi\'ee.

Avec cette approche, on peut aussi traiter le cas de deux (ou plusieurs) termes sources concentr\'es en $x_{1} \neq x_{2}$, deux points de $\R^{n}$. Pour $S^{1} , S^{2}$ deux fonctions satisfaisant l'hypoth\`ese \ref{h4}, on pose
\begin{equation*}
S_{h}^{1} (x) = h^{-n/2} S \Big( \frac{x- x_{1} }{h} \Big) \quad \text{et} \quad S_{h}^{2} (x) = h^{-n/2} S \Big( \frac{ x - x_{2}}{h} \Big) .
\end{equation*}
En plus de l'hypoth\`ese \ref{h5} aux points $x_{1}$ et $x_{2}$, on rajoute :
\begin{hyp}  \label{h8}
On suppose que
\begin{equation*}
\text{mes}_{n-1} \big\{ \xi \in \sqrt{2 (E_{0} - V(x_{1} ))} \S^{n-1} ; \ \exists t > 0 \quad \Pi_{x} \exp ( t {\rm H}_{p}) (x_{1} , \xi) = x_{2} \big\} =0 .
\end{equation*}
$\text{mes}_{n-1}$ est la mesure de Lebesgue sur $\sqrt{2 (E_{0} - V(x_{1}))} \S^{n-1}$.
\end{hyp}
Autrement dit on suppose que les trajectoires dans $p^{-1} ( E_{0})$ qui vont de $x_{1}$ \`a $x_{2}$ sont de mesure nulle. L'hypoth\`ese \ref{h8} ne change pas si $x_{1}$ et $x_{2}$ sont intervertis.

Pour $j =1 ,2$, soient $\Lambda_{k}^{j}$, $k = 0 , \ldots , K_{j}$ les vari\'et\'es lagrangiennes d\'efinies en \eqref{a79} et \eqref{a53} en $\rho_{0}$ avec point de base $x_{j}$. L'hypoth\`ese \ref{h8} est alors \'equivalente \`a
\begin{hypri} \label{h9}
Les intersections des vari\'et\'es $\Lambda^{1}_{k_{1}}$ et $\Lambda^{2}_{k_{2}}$ sont de mesure nulle :
\begin{equation*}
\text{mes}_{n} ( \Lambda_{k_{1}}^{1} \cap \Lambda_{k_{2}}^{2} ) =0 ,
\end{equation*}
pour tout $0 \leq k_{j} \leq K_{j}$, $j = 1,2$. La mesure est prise sur la vari\'et\'e $\Lambda_{k_{1}}^{1}$ (ou $\Lambda_{k_{2}}^{2}$).
\end{hypri}
Dans ce cas, la solution de Helmholtz est alors la fonction $u_{h}$ solution de
\begin{equation}
u_{h} (x) = (P- (E+ i 0))^{-1} \big( S_{h}^{1} + S_{h}^{2} \big) .
\end{equation}
On a alors le th\'eor\`eme suivant qui \'etend le r\'esultat de Fouassier \cite{Fo06_01} aux cas \`a coefficients variables.
\begin{theo}\sl
On suppose \ref{h1}, \ref{h2}, \ref{h6} et \ref{h8}. On suppose aussi que les hypoth\`eses \ref{h3}, \ref{h5} et \ref{h4} sont v\'erifi\'ees en $x_{1}$ et en $x_{2}$ (autrement dit, en rempla\c{c}ant le point de base $0$ par $x_{1}$ et par $x_{2}$). Alors, pour tout $q \in C^{\infty}_{0} (\R^{2n})$,
\begin{equation}
h \< \Op (q) u_{h} , u_{h} \> \longrightarrow \< q , \mu_{h^{-1}}^{1} +  \mu_{h^{-1}}^{2} \>  \quad \text{quand } h \to 0 .
\end{equation}
La mesure $\mu_{h^{-1}}^{j}$ est la mesure fournie par le th\'eor\`eme \ref{a48} avec point de base $x_{j}$.
\end{theo}

\noindent
La preuve suit celle du th\'eor\`eme \ref{a48} et on omet les d\'etails. L'hypoth\`ese \ref{h8} implique que les interactions entre les deux points sources sont n\'egligeables. Comme pour le th\'eor\`eme \ref{a48}, si cette hypoth\`ese est retir\'ee, on peut construire un exemple d'op\'erateur pour lequel la mesure n'est pas unique. Remarquons que ces hypoth\`eses impliquent $n \geq 2$.

Le reste du papier est consacr\'e \`a la preuve du th\'eor\`eme \ref{a48}. Dans la section \ref{a50}, on donne un exemple o\`u la mesure n'est pas unique. La proposition~\ref{a59} est d\'emontr\'ee dans la section \ref{a60}. Enfin, l'\'equivalence entre \ref{h5} et \ref{h7} est prouv\'ee dans la section \ref{a85}.

\section{Preuve du th\'eor\`eme \ref{a48}.}
\label{a63}

Dans les sections \ref{a62}--\ref{a64}, on pourra supposer que $E = E_{0} + h F_{1}$, o\`u $F_{1} \in \{ z \in \C ; \ \im z \geq 0 \}$ ne d\'epend pas de $h$. En effet, toutes les constructions seront r\'eguli\`eres par rapport au param\`etre $F_{1}$ et toutes les estimations seront uniformes en $F_{1}$, quand $F_{1}$ varie dans un compact de $\{ z \in \C ; \ \im z \geq 0 \}$. Et pour obtenir le cas g\'en\'eral \ref{h6}, il suffira de prendre $F_{1} = h^{-1} (E - E_{0}) = E_{1} + o(1)$, o\`u $E$ est donn\'e par \ref{h6}.

\Subsection{Localisation pr\`es de la surface d'\'energie $E_{0}$.}
\label{a62}

On note $P_{0} = - h^{2} \Delta /2$ le laplacien libre de symbole $p_{0} ( x , \xi ) = \xi^{2} /2$. 
\begin{lemm}\sl  \label{a4}
Soit $f,g \in C^{\infty}_{0} (\R)$ tel que $f=1$ pr\`es de $E_{0} - V(0)$ et $g =1$ pr\`es de $E_{0}$. Alors
\begin{align}
u_{h} (x) =& (P- (E+ i 0))^{-1} f (P_{0}) S_{h} + \CO^{\alpha} (1)  \label{a49} \\
=& g (P) (P- (E+ i 0))^{-1} f (P_{0}) S_{h} + \CO^{\alpha} (1) ,  \label{a45}
\end{align}
pour tout $\alpha >1/2$.
\end{lemm}

\begin{proof}
Soit $1_{E_{0}} \prec \widetilde{f} \in C^{\infty}_{0} ( \R)$ tel que $f (p_{0}) =1$ pr\`es du support de $\widetilde{f} (p)$. Par le calcul fonctionnel et comme $S_{h}= \CO (1)$, on a
\begin{align}
u_{h} (x) =& (P- (E+ i 0))^{-1} \widetilde{f} (P) S_{h} + (P - E)^{-1} (1 - \widetilde{f} (P)) S_{h}   \nonumber \\
=& (P- (E+ i 0))^{-1} \widetilde{f} (P) S_{h} + \CO (1) .  \label{a2}
\end{align}
Ici, on a utilis\'e le fait que $\widetilde{f} (P) \in \Psi ( 1)$ et donc $\widetilde{f} (P) S_{h} \in L^{2} ( \< x \>^{- \infty} dx)$. D'apr\`es le calcul fonctionnel des op\'erateurs pseudo-diff\'erentiels (\cite{DiSj99_01} par exemple), on a
\begin{equation}  \label{a3}
\widetilde{f} (P) \in \widetilde{f} (P) f (P_{0}) + \Psi (h^{\infty} \< \xi \>^{- \infty}) .
\end{equation}
En utilisant \eqref{a2}, \eqref{a3}, \eqref{a1} et encore \eqref{a2}, il vient
\begin{align}
u_{h} (x) =& (P- (E+ i 0))^{-1} \widetilde{f} (P) S_{h} + \CO (1)   \nonumber  \\
=& (P- (E+ i 0))^{-1} \widetilde{f} (P) f (P_{0}) S_{h} + \CO (1) + \CO^{\alpha} (h^{\infty}) \nonumber \\
=& (P- (E+ i 0))^{-1} f (P_{0}) S_{h} + \CO (1) + \CO^{\alpha} (h^{\infty}) ,
\end{align}
ce qui prouve \eqref{a49}. L'identit\'e \eqref{a45} d\'ecoule de \eqref{a2}.
\end{proof}

Remarquons que, comme $(P - E)^{-1} (1 - \widetilde{f} (P)) \in \Psi (\<\xi \>^{-2} )$, la preuve montre que
\begin{equation} \label{a77}
u_{h} (x) = (P- (E+ i 0))^{-1} f (P_{0}) S_{h} + R ,
\end{equation}
o\`u $R = \CO^{\alpha} (1)$ v\'erifie $(1 - \varphi (x)) R = \CO^{\alpha} (h^{\infty})$ pour toute fonction $\varphi \in C^{\infty}_{0} (\R^{n})$ telle que $\varphi =1$ pr\`es de $0$. En fait, comme $(P- (E+ i 0))^{-1} = \CO (h^{-1})$ en tant qu'op\'erateur de $H^{s} ( \< x \>^{2 \alpha} dx)$ dans $H^{s} ( \< x \>^{- 2 \alpha} dx)$, pour tout $s > 0$, on peut montrer que $\partial^{\beta}_{x} R = \CO^{\alpha} (h^{- \vert \beta \vert} )$ et $(1 - \varphi (x)) \partial^{\beta}_{x} R = \CO^{\alpha} (h^{\infty})$ pour tout $\beta \in \N^{n}$.

\Subsection{Contr\^ole des temps grands.}

Comme dans Robert--Tamura \cite{RoTa89_01} et Castella \cite{Ca05_01}, on voit la r\'esolvante comme l'int\'egrale en temps du propagateur. En effet,
\begin{equation}  \label{a5}
(P- (E+ i 0))^{-1} = \frac{i}{h} \int_{0}^{T} e^{-i t (P -E)/h} d t + (P- (E+ i 0))^{-1} e^{-i T (P - E)/h} ,
\end{equation}
pour tout $T\in \R$. On a alors le lemme suivant.

\begin{lemm}\sl  \label{a6}
Pour tout $q \in C^{\infty}_{0} (\R^{2 n})$,
\begin{equation}  \label{a18}
\Op ( q ) (P- (E+ i 0))^{-1} f (P_{0}) S_{h} = \Op ( q ) \frac{i}{h} \int_{0}^{T} e^{-i t (P -E)/h} f (P_{0}) S_{h} \, d t + \CO^{- \beta} (h^{\infty}) ,
\end{equation}
pour tout $T$ assez grand et $\beta > 0$.
\end{lemm}

\begin{proof}
La preuve est enti\`erement donn\'ee au bas de la page 173 de \cite{RoTa89_01}, et on rappelle juste les deux \'etapes. D'abord, on \'enonce une version du th\'eor\`eme d'Egorov.

\begin{prop}[Proposition 3.1 de \cite{RoTa89_01}]\sl  \label{a72}
Soit $\omega (x, \xi) \in C^{\infty}_{0} (\R^{2n})$. Supposons que $\omega_{t} \in S (1)$ soit nul pr\`es de
\begin{equation*}
\{ (x, \xi ) ; \ (x , \xi ) = \exp ( t {\rm H}_{p} ) ( y , \eta ) \text{ avec } (y , \eta ) \in \supp \omega \}.
\end{equation*}
Alors
\begin{equation*}
\big\Vert \< x \>^{\beta} \Op ( \omega_{t} ) e^{-i t P/h} \Op (\omega ) \< x \>^{\beta} \big\Vert = \CO (h^{\infty}) ,
\end{equation*}
pour tout $\beta \geq 0$. Cette estimation reste uniform\'ement vraie pour $t$ dans un compact de $\R$ \`a condition que $\omega_{t}$ soit uniform\'ement dans $S (1)$.
\end{prop}

\noindent
 On dit qu'un symbole $\omega \in C^{\infty} (\R^{2n})$ est de classe $A_{0}$ si
\begin{equation*}
\vert \partial_{x}^{\alpha} \partial_{\xi}^{\beta} \omega (x, \xi) \vert  \lesssim \< x \>^{- \vert \alpha \vert} \< \xi \>^{-L},
\end{equation*}
pour tout $L > 0$. Les ensembles $\Gamma_{-}$ (resp. $\Gamma_{+}$) d\'efinis par
\begin{equation*}
\Gamma_{\pm} ( R , d , \sigma ) = \{ (x , \xi ) ; \ \vert x \vert > R , \ d^{-1} < \vert \xi \vert < d \text{ et } \pm \cos ( x, \xi ) > \sigma \} ,
\end{equation*}
sont appel\'es les zones entrantes (resp. sortantes). Dans \cite{RoTa89_01}, la proposition \ref{a72} est \'enonc\'ee pour $\omega_{t} \in A_{0}$, mais elle reste vraie pour $\omega_{t} \in S (1)$.

Comme les surfaces d'\'energie proche de $E_{0}$ sont non-captives, cette proposition implique que le micro-support de $e^{-i T (P - E)/h} f (P_{0}) S_{h} (x)$ est micro-localis\'e dans une zone sortante modulo $\CO^{- \beta} (h^{\infty})$, pour tout $\beta \geq 0$. Le lemme \ref{a6} est alors une cons\'equence directe du lemme suivant qui repose sur les constructions d'Isozaki et Kitada.

\begin{prop}[Lemme 2.3 de \cite{RoTa89_01}]\sl  \label{a43}
Soit $\omega$, $\omega_{\pm} \in A_{0}$ tel que le support de $\omega_{\pm}$ (resp. $\omega$) soit dans $\Gamma_{\pm} (R , d , \sigma_{\pm})$ (resp. $\{ \vert x \vert < 9 R /10 \}$), avec $R$ assez grand et $\sigma_{-} < \sigma_{+}$. Alors
\begin{gather*}
\big\Vert \< x \>^{\beta} \Op ( \omega_{\mp} ) (P - (E \pm i 0))^{-1} \Op (\omega_{\pm} ) \< x \>^{\beta} \big\Vert = \CO (h^{\infty}) , \\
\big\Vert \< x \>^{\beta} \Op ( \omega ) (P - (E \pm i 0))^{-1} \Op (\omega_{\pm} ) \< x \>^{\beta} \big\Vert = \CO (h^{\infty}) ,
\end{gather*}
pour tout $\beta \geq 0$.
\end{prop}

\noindent
Remarquons que Robert et Tamura suppose que l'\'energie $E$ est r\'eelle et que le potentiel est \`a courte port\'ee, $\rho >1$ (pour pouvoir d\'evelopper la th\'eorie de la diffusion), mais les propositions pr\'ec\'edentes sont encore valables pour des potentiels \`a longue port\'ee (voir l'appendice de \cite{RoTa88_01} o\`u la construction des propagateurs approch\'es entrant et sortant est faites pour $\rho > 0$ (voir aussi les travaux d'Isozaki et Kitada \cite{IsKi85_01} et \cite{IsKi85_02})) et pour des \'energies v\'erifiant \ref{h6} (en fait, pour $E \in \C$ proche de $E_{0}$ et $\pm \im E \geq 0$).
\end{proof}

\Subsection{Contr\^ole des temps petits.} \label{a75}

En utilisant les lemmes \ref{a4} et \ref{a6} et les propri\'et\'es du calcul pseudo-diff\'erentiel, on obtient, pour $q \in C^{\infty}_{0} ( \R^{2n} )$,
\begin{equation}  \label{a76}
\Op (q) u_{h} = \Op (q) \widetilde{u}_{h} + \CO^{- \beta} (1),
\end{equation}
o\`u $\widetilde{u}_{h}$ est d\'efini par
\begin{equation}
\widetilde{u}_{h} := \frac{i}{h} \int_{0}^{T} e^{-i t (P -E)/h} f(P_{0}) S_{h} \, d t .
\end{equation}
Dans cette section, on \'etudie
\begin{equation}
B_{1} (x) := \frac{i}{h} \int_{0}^{T} \chi (t) e^{-i t (P -E)/h} f (P_{0}) S_{h} \, d t ,
\end{equation}
avec $\chi \in C^{\infty}_{0} (\R)$.

Remarquons que, pour toute fonction $g \in C^{\infty}_{0} (\R^{n})$ telle que $g =1$ pr\`es de $0$, 
\begin{align}
f(P_{0}) S_{h} =& g (x) f(P_{0}) S_{h} + \CO (h^{\infty})   \nonumber  \\
=& (2 \pi )^{-n} h^{- \frac{n}{2}} \int g (x) e^{i x \cdot \xi / h}  f ( \xi^{2} / 2) \widehat{S} ( \xi) \, d \xi + \CO (h^{\infty}) .  \label{a26}
\end{align}
De plus, le support de $g$ (resp. $f$) est aussi proche de $0$ (resp. $E_{0} - V(0)$) que l'on veut. Donc
\begin{equation}
B_{1} (x) = \frac{i h^{-1-\frac{n}{2}}}{(2 \pi)^{n}} \iint_{0}^{T} \chi (t) e^{-i t (P -E)/h} e^{i x \cdot \xi / h} g (x) f ( \xi^{2} / 2) \widehat{S} ( \xi) \, d t \, d \xi  + \CO (h^{\infty}) .
\end{equation}
La norme $L^{2}$ de $B_{1}$ tr\`es pr\`es de $x = 0$ est contr\^ol\'ee par le lemme suivant.

\begin{lemm}\sl \label{a7}
Soit $\chi \in C^{\infty}_{0} (\R)$, $\chi =1$ pr\`es de $0$ et de support suffisamment proche de $0$. Alors pour tout $\varepsilon >0$,
\begin{equation}
1_{\vert x \vert < \varepsilon} B_{1} = \CO \big( \sqrt{\varepsilon} h^{- \frac{1}{2}} \big) ,
\end{equation}
en norme $L^{2}$.
\end{lemm}

\begin{proof}
D'apr\`es la m\'ethode BKW (voir \cite{FeMa81_01} ou \cite{DiSj99_01} par exemple), on sait que, pour tout $t$ assez petit,
\begin{equation}
e^{-i t (P -E)/h} \big( g ( y ) e^{i y \cdot \xi / h} \big) = a (t, x, \xi ,h) e^{i \varphi (t, x, \xi) /h} + \CO (h^{\infty}) ,
\end{equation}
La phase $\varphi$ et le symbole $a$ sont $C^{\infty}$ par rapport \`a $t,x, \xi$. La phase v\'erifie l'\'equation eikonale
\begin{equation}  \label{a8}
\left\{ \begin{aligned}
&\partial_{t} \varphi + p (x, \partial_{x} \varphi ) - E_{0} = 0  \\
&\varphi (0, x, \xi ) = x \cdot \xi
\end{aligned} \right.
\end{equation}
L'amplitude $a (t,x, \xi)$ est \`a support compact en $x$ et d\'eveloppable par rapport \`a $h$
\begin{equation*}
a (t, x, \xi ,h) \sim \sum_{j =0}^{\infty} h^{j} a_{j} (t ,x , \xi).
\end{equation*}
De plus, les coefficients $a_{j}$ v\'erifient les \'equations de transport habituelles
\begin{equation}
\begin{aligned} \label{a47}
(\partial_{t} + \nabla_{x} \varphi \cdot \nabla_{x} + \Delta_{x} \varphi /2  - F_{1} ) a_{0} =& 0   \\
(\partial_{t} + \nabla_{x} \varphi \cdot \nabla_{x} + \Delta_{x} \varphi /2  - F_{1} ) a_{j} =&  \frac{1}{2} \Delta a_{j - 1} \quad \text{pour } j \geq 1 ,
\end{aligned}
\end{equation}
avec donn\'ees en $t=0$
\begin{equation}
\begin{aligned} \label{a87}
a_{0} (0 , x) =& g (x)  \\
a_{j} (0 ,x) =&0 \quad \text{pour } j \geq 1 .
\end{aligned}
\end{equation}
En particulier, comme annonc\'e au d\'ebut de la partie \ref{a63}, les $a_{j}$ sont holomorphes par rapport \`a $F_{1}$.

On est alors amen\'e \`a calculer
\begin{equation}
J := \iint_{0}^{T} \chi (t) e^{i \varphi (t, x, \xi) /h} a (t, x, \xi ,h) f ( \xi^{2} / 2) \widehat{S} ( \xi ) \, d \xi \, d t ,
\end{equation}
puisque $B_{1} = \frac{i h^{-1 - \frac{n}{2}}}{(2 \pi )^{n}} J + \CO (h^{\infty})$. Pour $\delta \gg h$ petit, on note
\begin{equation}
J_{\delta} (x) := 1_{\delta < \vert x \vert < 2 \delta} \iint_{0}^{T} \chi (t) e^{i \varphi (t, x, \xi) /h} a (t, x, \xi ,h) f ( \xi^{2} / 2) \widehat{S} ( \xi ) \, d \xi \, d t ,
\end{equation}
et on essaye d'estimer $J_{\delta}$.

%
%

Bien sûr, pour $\widetilde{\chi} \in C^{\infty}_{0} (\R)$,
\begin{equation}
J_{\delta}^{\infty} (x) := 1_{\vert x \vert < \delta} \iint_{0}^{T} \widetilde{\chi} (t / \delta ) \chi (t) e^{i \varphi (t, x, \xi) /h} a (t, x, \xi ,h) f ( \xi^{2} / 2) \widehat{S} ( \xi ) \, d \xi \, d t ,
\end{equation}
v\'erifie
\begin{equation}  \label{a14}
\Vert J_{\delta}^{\infty} \Vert_{L^{2}} = \CO ( \delta^{1 + \frac{n}{2}} ).
\end{equation}

%

Remarquons que, d'apr\`es \eqref{a8}, on a
\begin{equation} \label{a9}
\varphi (t, x, \xi) = x \cdot \xi - t ( \xi^{2} /2 + V (x) -E_{0} ) + t^{2} S (1) .
\end{equation}
De plus, grâce au lemme \ref{a4}, $\xi^{2}$ reste dans un compact disjoint de $0$ sur le support de $f (\xi^{2} / 2)$. En particulier, pour $C_{1} \delta < t < 1/ C_{1}$ avec $C_{1}$ assez grand (ind\'ependant de $\delta$), on~a
\begin{equation}
- \xi \cdot \partial_{\xi} \varphi = - x \cdot \xi + t \xi^{2} + t^{2} S(1) \geq \delta ,
\end{equation}
sur le support de $1_{\delta < \vert x \vert < 2 \delta} a (t,x, \xi ,h) f(\xi^{2}/2)$. Comme les d\'eriv\'ees par rapport \`a $\xi$ de la fonction $\frac{- \delta \xi}{-x \cdot \xi + t \xi^{2} + t^{2} S(1)}$ sont uniform\'ement born\'ees, des int\'egrations par partie en $\xi$ donnent, pour $\widetilde{\chi} \in C^{\infty}_{0} (] - \infty , C_{1} +1 ])$ avec $\widetilde{\chi} =1$ pr\`es de $] - \infty , C_{1} ]$,
\begin{align}
J_{\delta}^{+} (x) :=& 1_{\delta < \vert x \vert < 2 \delta} \iint_{0}^{T} \big( 1 - \widetilde{\chi} (t / \delta ) \big) \chi (t) e^{i \varphi (t, x, \xi) /h} a (t, x, \xi ,h) f ( \xi^{2} / 2) \widehat{S} ( \xi ) \, d \xi \, d t    \nonumber  \\
=& \CO \big( (h / \delta )^{M} \big),
\end{align}
pour tout $M \in \N$. Ainsi
\begin{equation}  \label{a15}
\Vert J_{\delta}^{+} \Vert_{L^{2}} = \CO \big( \delta^{\frac{n}{2}} (h / \delta )^{M} \big) ,
\end{equation}
pour tout $M \in \N$.

%

D'un autre cot\'e, pour $t < \delta / C_{2}$ avec $C_{2}$ assez grand, on a
\begin{equation}
\frac{x \cdot \partial_{\xi} \varphi}{x^{2}} = \frac{x^{2} - t x \cdot \xi + t^{2} x S (1)}{x^{2}} \geq 1 ,
\end{equation}
sur le support de $1_{\delta < \vert x \vert < 2 \delta} a (t,x, \xi ,h) f(\xi^{2}/2)$. Comme les d\'eriv\'ees par rapport \`a $\xi$ de la fonction $\frac{\delta x}{x^{2} - t x \cdot \xi + t^{2} x S (1)}$ sont uniform\'ement born\'ees, des int\'egrations par partie en $\xi$ donnent, pour $\overline{\chi} \in C^{\infty}_{0} (] - \infty , C_{2} ])$ avec $\overline{\chi} =1$ pr\`es de $] - \infty , C_{2}/2 ]$,
\begin{align}
J_{\delta}^{-} (x) :=& 1_{\delta < \vert x \vert < 2 \delta} \iint_{0}^{T} \overline{\chi} (t / \delta ) \chi (t) e^{i \varphi (t, x, \xi) /h} a (t, x, \xi ,h) f ( \xi^{2} / 2) \widehat{S} ( \xi ) \, d \xi \, d t    \nonumber  \\
=& \big( (h / \delta )^{M} \big) ,
\end{align}
pour tout $M \in \N$. Ainsi
\begin{equation}  \label{a16}
\Vert J_{\delta}^{-} \Vert_{L^{2}} = \CO \big( \delta^{\frac{n}{2}} (h / \delta )^{M} \big) ,
\end{equation}
pour tout $M \in \N$.

%

Il reste \`a \'etudier
\begin{equation}  \label{a19}
J_{\delta}^{0} (x) := 1_{\delta < \vert x \vert < 2 \delta} \iint_{0}^{T} \breve{\chi} (t / \delta ) \chi (t) e^{i \varphi (t, x, \xi) /h} a (t, x, \xi ,h) f ( \xi^{2} / 2) \widehat{S} ( \xi ) \, d \xi \, d t ,
\end{equation}
o\`u $\breve{\chi} = \big( 1 - \overline{\chi} \big) \widetilde{\chi}$ \`a son support compact disjoint de $0$. On fait alors le changement de variables
\begin{equation*}
t = \delta T , \qquad
x = \delta X , \qquad
\xi = \Xi .
\end{equation*}
Le symbole $1_{\delta < \vert x \vert < 2 \delta} \breve{\chi} (t / \delta ) \chi (t) a (t, x, \xi ,h) f ( \xi^{2} / 2) \widehat{S} ( \xi )$ devient $1_{1 < \vert X \vert < 2} \breve{a} ( T, X , \Xi )$ o\`u le symbole $\breve{a} \in S (1)$ est \`a support compact dans $]0, + \infty [_{T} \times ( \R^{n} \setminus \{ 0 \})_{X} \times ( \R^{n} \setminus \{ 0 \})_{\Xi}$. On pose
\begin{equation}
\Phi (T , X , \Xi , \delta ) = \delta^{-1}  \varphi (\delta T , \delta  X , \Xi) .
\end{equation}
Cette fonction est $C^{\infty}$ par rapport aux variables $T , X , \Xi$ \underline{et} $\delta$. D'apr\`es \eqref{a9}, elle s'\'ecrit
\begin{equation}  \label{a10}
\Phi (T , X , \Xi , \delta ) = X \cdot \Xi - T ( \Xi^{2}/2 + V(0) -E_{0} ) + \delta R (T , X , \Xi , \delta ) ,
\end{equation}
o\`u $R \in S (1)$ sur le support de $\breve{a}$. Toujours sur le support de $\breve{a}$, les points critiques en $(T, \Xi)$ de $\Phi$ sont les points tels que
\begin{equation}
\partial_{T , \Xi } \Phi = \left( \begin{aligned}
&- \Xi^{2}/2 - V(0) + E_{0} + \delta \partial_{T} R    \\
&X - T \Xi + \delta \partial_{\Xi} R   
\end{aligned} \right) = 0 .
\end{equation}
Pour $\delta =0$, il existe un unique point critique $\Big( \widetilde{T} (X ,0) = \frac{\vert X \vert}{\sqrt{2 (E_{0} - V (0))}}, \widetilde{\Xi} (X ,0) = \frac{X \sqrt{2(E_{0} - V(0))}}{\vert  X \vert} \Big)$ (car $T >0$ sur $\supp (\breve{a})$). En ce point, 
\begin{equation}  \label{a11}
\Hess_{T , \Xi } \Phi = \left( \begin{array}{cc}
0 & -^{t} \Xi \\
- \Xi & -T
\end{array} \right) ,
\end{equation}
est de d\'eterminant $- T^{n-1} \Xi^{2}$ ce qui est uniform\'ement minor\'e sur $\supp ( \breve{a})$. Par le th\'eor\`eme des fonctions implicite, pour $\delta$ assez petit, l'\'equation $\partial_{T , \Xi } \Phi =0$ a donc une unique solution $( \widetilde{T}, \widetilde{\Xi} ) (X, \delta) = ( \widetilde{T}, \widetilde{\Xi} ) (X, 0) + \delta \widetilde{R}$ avec $\widetilde{R} \in S (1)$. En particulier, d'apr\`es \eqref{a10} et \eqref{a11}, $\Hess_{T , \Xi } \Phi$ est uniform\'ement minor\'e au point critique pour $\delta \ll 1$. En \'ecrivant
\begin{equation}
\partial_{T , \Xi } \Phi ( T , \Xi ) = \Hess_{T , \Xi } \Phi ( \widetilde{T} , \widetilde{\Xi} ) ( T - \widetilde{T} , \Xi - \widetilde{\Xi}) +  \CO ( \vert T - \widetilde{T} \vert^{2} , \vert \Xi - \widetilde{\Xi} \vert^{2}) ,
\end{equation}
on constate que
\begin{equation}
\frac{\vert (T , \Xi ) - (\widetilde{T} , \widetilde{\Xi}) \vert}{\vert \partial_{T , \Xi } \Phi (T , \Xi) \vert} ,
\end{equation}
est uniform\'ement born\'e. En appliquant la m\'ethode de la phase stationnaire (le param\`etre semi-classique est ici $h/\delta$) \`a param\`etre (th\'eor\`eme 7.7.5 de \cite{Ho90_01}), il vient
\begin{align}
J_{\delta}^{0} (x) =& \delta \iint 1_{1 < \vert X \vert < 2} \breve{a} ( T, X , \Xi ) e^{i \Phi (T , X , \Xi , \delta ) \delta /h} d \Xi \, d T   \nonumber  \\
=& \CO \big( \delta (h / \delta)^{\frac{n+1}{2}} \big) .   \label{a12}
\end{align}
Et donc
\begin{equation} \label{a13}
\Vert J_{\delta}^{0} (x) \Vert_{L^{2}} = \CO \big( \delta^{\frac{n}{2}} \delta (h / \delta)^{\frac{n+1}{2}} \big) = \CO \big( \delta^{\frac{1}{2}} h^{\frac{n+1}{2}} \big) .
\end{equation}

En, utilisant une d\'ecomposition dyadique $\delta = 2^{-n}$ avec $h^{1- \mu } < \delta < 2 \varepsilon$ ainsi que \eqref{a14}, \eqref{a15}, \eqref{a16} et \eqref{a13}, il vient
\begin{align}
\Vert 1_{\vert x \vert < \varepsilon} J \Vert_{L^{2}} \leq& \big( \Vert J_{h^{1- \mu}}^{\infty} \Vert + \Vert J_{h^{1- \mu}}^{+} \Vert \big) + \Big( \sum_{\fract{\delta = 2^{-n}}{h^{1- \mu } < \delta < 2 \varepsilon}} \Vert J_{\delta}^{-} \Vert + \Vert J_{\delta}^{0} \Vert + \Vert J_{\delta}^{+} \Vert \Big)   \nonumber  \\
\lesssim& h^{(1- \mu)(\frac{n}{2} +1)} + h^{(1- \mu)\frac{n}{2}} h^{\mu M} + \sum_{\fract{\delta = 2^{-n}}{h^{1- \mu } < \delta < 2 \varepsilon}} \delta^{\frac{n}{2}} (h / \delta )^{M} + \delta^{\frac{1}{2}} h^{\frac{n+1}{2}} + \delta^{\frac{n}{2}} (h / \delta )^{M}   \nonumber   \\
\lesssim& h^{\frac{n+1}{2}} \big( h^{\frac{1}{2} - \mu (\frac{n}{2} +1)} + h^{\mu M -\frac{1}{2} - \mu  \frac{n}{2}} + h^{\mu M - \frac{n+1}{2}} \vert \ln h \vert + \sqrt{\varepsilon} + h^{\mu M - \frac{n+1}{2}} \vert \ln h \vert)  \nonumber  \\
\lesssim& h^{\frac{n+1}{2}} ( \sqrt{\varepsilon} + h^{\nu}) ,  \label{a17}
\end{align}
avec $\nu >0$. Pour obtenir \eqref{a17}, on fixe $\mu >0$ assez petit (pour majorer le premier terme) puis $M$ assez grand.
\end{proof}

\Subsection{Calcul pr\`es de $0$.} \label{a74}

Dans cette partie, on reprend le calcul de
\begin{equation}
B_{1} (x) = \frac{i h^{-1-\frac{n}{2}}}{(2 \pi)^{n}} \iint_{0}^{T} \chi (t) e^{i \varphi (t, x, \xi) /h} a (t, x, \xi ,h) f ( \xi^{2} / 2) \widehat{S} ( \xi ) \, d \xi \, d t + \CO ( h^{\infty} ) ,
\end{equation}
pour $\delta_{1} < \vert x \vert < \delta_{2}$ avec $0< \delta_{1} < \delta_{2}$ petits mais fix\'es. D'apr\`es \eqref{a16}, on a vu que
\begin{equation}  \label{a20}
B_{1} (x) = \frac{i h^{-1-\frac{n}{2}}}{(2 \pi)^{n}} \iint_{0}^{T} \breve{\chi} (t) e^{i \varphi (t, x, \xi) /h} a (t, x, \xi ,h) f ( \xi^{2} / 2) \widehat{S} ( \xi ) \, d \xi \, d t + \CO ( h^{\infty} ) ,
\end{equation}
avec $\breve{\chi} = \chi ( 1 - \overline{\chi})$, o\`u $\overline{\chi} \in C^{\infty}_{0} (\R)$ et $\overline{\chi}=1$ pr\`es de $0$. D'apr\`es \eqref{a19}--\eqref{a12}, on sait que l'on peut utiliser la m\'ethode de la phase stationnaire en $t, \xi$ dans \eqref{a20} (pour tout $x$, la phase a un unique point critique qui est non d\'eg\'en\'er\'e). D'un autre cot\'e, en un tel point, on a
\begin{equation}  \label{a21}
\left\{ \begin{aligned}
&\partial_{t} \varphi (t, x, \xi ) =0  \\
&\partial_{\xi} \varphi (t, x, \xi ) =0 .
\end{aligned} \right.
\end{equation}
Rappelons que, d'apr\`es la proposition IV-14 i) de \cite{Ro87_01},
\begin{equation}   \label{a23}
( x, \partial_{x} \varphi ) = \exp (t {\rm H}_{p} ) (\partial_{\xi} \varphi , \xi ) ,
\end{equation}
et d'apr\`es (40) de \cite{Ro87_01}, on a
\begin{equation}  \label{a25}
\xi (s) = \partial_{x} \varphi (s, x (s) ,\xi ) ,
\end{equation}
pour toute courbe hamiltonienne $(x (s) , \xi (s))$ avec $\xi (0) = \xi$. En utilisant \eqref{a8}, \eqref{a23} et \eqref{a25}, on voit que \eqref{a21} est \'equivalente \`a
\begin{equation}  \label{a22}
\left\{ \begin{aligned}
&p (0, \xi ) = E_{0}  \\
&\Pi_{x} \exp (t {\rm H}_{p} ) (0 , \xi ) =x .
\end{aligned} \right.
\end{equation}
En utilisant \eqref{a25} et en notant $\gamma (x) = ( x(s) , \xi (s) )$ la courbe hamiltonienne $\exp (s {\rm H}_{p}) (0 , \xi)$ donn\'ee par \eqref{a22}, il vient
\begin{align}
\psi (x) :=& \varphi (t, x, \xi ) = \varphi (0 , x (0) , \xi ) + \int_{0}^{t} \partial_{x} \varphi (s, x (s) , \xi) \dot{x} (s) \, d s  \nonumber \\
=& \int_{0}^{t} \xi (s) \dot{x} (s) \, d s = \int_{\gamma (x)} \xi \, d x.
\end{align}
Et donc, la m\'ethode de la phase stationnaire dans \eqref{a20} donne
\begin{equation}  \label{a24}
B_{1} (x) = b (x,h) h^{-\frac{1}{2}} e^{i \psi (x) /h} + \CO (h^{\infty}) ,
\end{equation}
o\`u $b (x,h) \sim \sum_{j \geq 0} b_{j}(x) h^{j}$ est un symbole classique \`a support compact avec
\begin{equation}  \label{a31}
b_{0} = i  (2 \pi )^{\frac{1-n}{2}} \frac{e^{i \frac{\pi}{4} \sgn \varphi''_{t, \xi}}}{\vert \det \varphi''_{t, \xi} \vert^{\frac{1}{2}}} \breve{\chi} (t) a_{0} (t, x, \xi ) f ( \xi^{2} / 2) \widehat{S} ( \xi ) ,
\end{equation}
o\`u $(t, \xi )$ est l'unique couple de points (avec $t$ dans le support de $\breve{\chi}$) v\'erifiant \eqref{a22}. A l'aide de \eqref{a11}, on trouve $\sgn \varphi''_{t, \xi} = 1-n$ et on remarque que
\begin{equation}  \label{a32}
\vert \det \varphi''_{t, \xi} \vert^{\frac{1}{2}} = t^{\frac{n-1}{2}} \vert \xi \vert (1 + o(1)) ,
\end{equation}
quand $t \to 0$.

Notons $(x( t, \xi ) , \xi (t, \xi ))$ la courbe hamiltonienne issue de $(0 , \xi)$ avec $p (0, \xi ) =E_{0}$. La d\'efinition du champ hamiltonien donne $x (t, \xi ) = t \xi + t^{2} S (1)$ et, d'apr\`es la discussion pr\'ec\'edente, $(t, \xi ) \to x (t, \xi )$ est un diff\'eomorphisme pour $t$ assez petit mais disjoint de $0$. De plus, son jacobien v\'erifie :
\begin{equation} \label{a33}
\left\vert \frac{d x}{ d (t, \xi )} \right\vert = \xi^{2} t^{n-1} (1 + \CO (t)) .
\end{equation}

\Subsection{Calcul pour les temps interm\'ediaires.}
\label{a55}

Dans cette partie, on calcule
\begin{equation}
\widetilde{u}_{h} = \frac{i}{h} \int_{0}^{T} e^{-i t (P -E)/h} f(P_{0}) S_{h} \, d t ,
\end{equation}
micro-localement pr\`es de $\rho_{0} \in p^{-1} (E_{0})$. Autrement dit, on calcule $\Op (q) \widetilde{u}_{h}$, pour $q$ un symbole \`a support compact proche de $\rho_{0}$.

Soient $t_{1} , \ldots , t_{K}$ les temps strictement positifs d\'efinies dans \eqref{a69}. On se donne une fonction $\chi \in C^{\infty}_{0} (\R)$ \'egale \`a $1$ pr\`es de $0$ et \`a support assez proche de $0$ (et donc $\chi$ est nulle pr\`es des $t_{k}$). Soit $\breve{\chi} \in C^{\infty}_{0} (] 0 , + \infty [)$ une fonction comme dans la partie pr\'ec\'edente, \`a support assez proche de $0$ et v\'erifiant $\breve{\chi} ( t -\nu ) = 1$ pr\`es de $0$ pour un certain $\nu >0$ petit. Il existe alors une fonction $\widetilde{\chi} \in C^{\infty}_{0} (\R )$, nulle pr\`es de $0$ et des $t_{k}$, telle que, pour $t \in [ 0 , T]$,
\begin{equation}
\chi (t) + \sum_{k =1}^{K} \breve{\chi} ( t -\nu - t_{k} ) + \widetilde{\chi} (t) = 1.
\end{equation}
Ainsi, $\widetilde{u}_{h}$ s'\'ecrit
\begin{equation}  \label{a73}
\widetilde{u}_{h} = B_{1} + B_{2} + \frac{i}{h} \int_{0}^{T} \widetilde{\chi} (t) e^{-i t (P -E)/h} f(P_{0}) S_{h} \, d t ,
\end{equation}
avec
\begin{equation}
B_{2} := \sum_{k = 1}^{K} \frac{i}{h} \int \breve{\chi} (t - \nu - t_{k} ) e^{-i t (P -E)/h} f(P_{0}) S_{h} \, d t .
\end{equation}

Rappelons que, d'apr\`es \cite{Ro87_01}, l'\'evolution est un op\'erateur int\'egral de Fourier d'ordre $0$ :
\begin{equation}  \label{a27}
e^{-i t (P-E) /h} \in \CI^{0}_{{\rm cl}} ( \Lambda_{t}) ,
\end{equation}
dont la relation canonique est donn\'ee par le flot hamiltonien :
\begin{equation}   \label{a28}
\begin{gathered}
\Lambda_{t} = \{ (x, \xi , y , \eta ) ; \ (x, \xi ) = \kappa_{t} (y, \eta ) \}  \\
\kappa_{t} ( y , \eta ) = \exp (t {\rm H}_{p} ) (y, \eta) .
\end{gathered}
\end{equation}
La notation ${\rm cl}$ signifie que le symbole de l'OIF \`a un d\'eveloppement en puissance de $h$. De plus, le symbole et chaque terme de son d\'eveloppement sont holomorphes par rapport \`a $F_{1}$.

Comme le micro-support de $f (P_{0} ) S_{h}$ est contenu dans $\{ (0 , \xi) ; \ \xi^{2} /2 \in \supp f \}$, le th\'eor\`eme d'Egorov implique que le dernier terme dans \eqref{a73} est micro-localement nul pr\`es de $\rho_{0}$. En effet, pour $t \in \supp {\bf 1}_{[0 ,T]} \widetilde{\chi} (t)$, l'op\'erateur d\'evolution $e^{-i t (P -E)/h}$ envoie le micro-support de $f(P_{0}) S_{h}$ hors du support de $q$.  Rappelons que ``$A = B$ micro-localement pr\`es de $\rho_{0}$'' signifie que $\Op (q) A = \Op (q) B + \CO (h^{\infty} )$. Comme $B_{1}$ est bien d\'ecrit par les sections \ref{a75} et \ref{a74}, il nous reste \`a \'etudier $B_{2}$.

On peut \'ecrire
\begin{equation}
B_{2} = \sum_{k =0}^{K} e^{-i ( \nu + t_{k} ) (P -E)/h} \Big( \frac{i}{h} \int_{0}^{+ \infty} \breve{\chi}(t) e^{-i t (P -E)/h} f(P_{0}) S_{h} \, d t \Big) .  \label{a70}
\end{equation}
D'apr\`es la partie pr\'ec\'edente, on sait que
\begin{equation}
\frac{i}{h} \int_{0}^{+ \infty} \breve{\chi} (t) e^{-i t (P -E)/h} f(P_{0}) S_{h} \, d t = h^{-\frac{1}{2}} b (x ,h) e^{i \psi (x) /h} + \CO ( h^{\infty} ) ,
\end{equation}
o\`u $b (x,h)$ est un symbole classique \`a support dans un compact proche de $0$. En particulier,
\begin{equation} \label{a29}
b (x ,h) e^{i \psi (x) /h} \in \dl^{0}_{\text{cl}} ( \widetilde{\Lambda} , b_{0} , \Omega ),
\end{equation}
o\`u $\Omega \subset \R^{2n}$ est un petit voisinage de $\{ (0 , \xi ) ; \ \xi^{2} = 2 (E_{0} - V (0)) \}$. La notation $\dl^{\alpha}_{\text{cl}} ( \Lambda , a_{0} , U )$ d\'esigne l'ensemble des distributions lagrangiennes semi-classiques, $u$, d'ordre $\alpha$, de vari\'et\'e lagrangienne $\Lambda$, de symbole principal $a_{0}$ (d\'efini sur le fibr\'e de Maslov), de symbole classique (c'est \`a dire que dans toute repr\'esentation \`a l'aide d'une phase non-d\'eg\'en\'er\'ee, le symbole de l'int\'egrale oscillante $a (x, \theta ,h)$ v\'erifie $a \sim \sum_{j \geq 0} a_{j}(x, \theta ) h^{j}$) et de micro-support inclus dans $U \subset \R^{2 n}$ (c'est \`a dire $\Op ( 1 - q ) u = \CO (h^{\infty} )$ pour tout $q \in S (1)$ avec $q =1$ pr\`es de $U$). Pour toute les questions concernant les distributions lagrangiennes, on renvoie \`a la section 1.2.1 de \cite{Iv98_01} pour un expos\'e semi-classique et \`a \cite{Ho94_01} dans le cas classique.

Mais alors, d'apr\`es \eqref{a27}, \eqref{a28}, \eqref{a29}, le th\'eor\`eme d'Egorov et le th\'eor\`eme 1.2.12 {\sl (iii)} de \cite{Iv98_01} sur le calcul de l'application d'un OIF \`a une int\'egrale oscillante, on a
\begin{equation*}
e^{-i ( \nu + t_{k} )(P -E)/h} h^{-\frac{1}{2}} b (x ,h) e^{i \psi (x) /h} \in h^{-\frac{1}{2}} \dl^{0}_{\text{cl}} \big( \exp (( \nu + t_{k} ) {\rm H}_{p}) \widetilde{\Lambda} , b^{k}_{0} , \exp (( \nu + t_{k} ) {\rm H}_{p}) (\Omega) \big) .
\end{equation*}
En notant $\Lambda_{k} := \exp ( ( \nu + t_{k} ) {\rm H}_{p}) \widetilde{\Lambda} \subset \Lambda$ et $\Omega_{k} := \exp ( ( \nu + t_{k} ) {\rm H}_{p}) (\Omega)$, on a donc,
\begin{equation}  \label{a54}
B_{2} = \sum_{k=0}^{K} B_{2}^{k} \quad \text{avec} \quad  B_{2}^{k} \in h^{-\frac{1}{2}} \dl^{0}_{\text{cl}} ( \Lambda_{k} , b^{k}_{0} , \Omega_{k} ) ,
\end{equation}
micro-localement pr\`es de $\rho_{0}$. Les $\Lambda_{k}$ coïncident avec ceux d\'efinis en \eqref{a53} et les symboles des distributions lagrangiennes ainsi que les $b^{k}_{0}$ sont r\'eguliers par rapport \`a $F_{1}$.

Finalement remarquons que, micro-localement pr\`es de $\rho_{0}$, on a
\begin{align}
(P-E) B_{2}^{k} =& (P-E) \frac{i}{h} \int_{0}^{+ \infty} \breve{\chi} (t - \nu - t_{k} ) e^{-i t (P -E)/h} f(P_{0}) S_{h} \, d t  \nonumber  \\
=& - \int_{\R} \breve{\chi} (t - \nu - t_{k} ) \partial_{t} e^{-i t (P -E)/h} f(P_{0}) S_{h} \, d t  \nonumber  \\
=& \int_{\R} \partial_{t} (\breve{\chi}) (t - \nu - t_{k} ) e^{-i t (P -E)/h} f(P_{0}) S_{h} \, d t .
\end{align}
En utilisant le th\'eor\`eme d'Egorov et les propri\'et\'es du support de $\partial_{t} (\breve{\chi})$, il vient
\begin{equation} \label{a71}
(P-E) B_{2}^{k} = 0 ,
\end{equation}
micro-localement pr\`es de $\rho_{0}$.

\Subsection{Convergence vers une mesure de Wigner.}
\label{a64}

\begin{prop}\sl  \label{a39}
Il existe une mesure de Radon positive $\mu_{h^{-1}}$ sur $\R^{2n}$ telle que, pour tout $q \in C^{\infty}_{0} (\R^{2n})$,
\begin{equation}
h \< \Op (q) u_{h} , u_{h} \> \longrightarrow \< q , \mu_{h^{-1}} \> ,
\end{equation}
quand $h \to 0$.
\end{prop}

\begin{proof}
Par lin\'earit\'e, on peut supposer que $q$ est \`a support compact pr\`es d'un point $\rho_{0}$. Soit $v (x) \in C^{\infty}_{0} ( \R^{n} )$ avec $v =1$ pr\`es de $0$ et $v_{\delta} (x) = v (x/ \delta )$ pour $\delta > 0$. D'apr\`es l'hypoth\`ese \ref{h7} et comme l'intersection de deux vari\'et\'es $\Lambda_{j}$ et $\Lambda_{k}$ est forc\'ement un ferm\'e, la r\'egularit\'e de la mesure combin\'ee avec le lemme d'Urysohn implique qu'il existe une suite de fonctions $f_{m} (x, \xi ) \in C^{\infty}_{0} (\R^{2n} , [0,1])$ telles que, sur le support de $q$,
\begin{gather}
\text{mes}_{n} ( \supp (1-f_{m}) \cap \Lambda ) \longrightarrow 0 ,
\end{gather}
lorsque $m \to + \infty$ et $\Lambda$ restreint au support de $f_{m}$ ne se recoupe pas.

Plus explicitement, les $f_{m}$ sont construites de la fa\c{c}on suivante : comme le ferm\'e $\Lambda_{j} \cap \Lambda_{k}$ est de mesure nulle dans la vari\'et\'e lisse $\Lambda_{j}$, la r\'egularit\'e de la mesure entraine l'existence d'un ouvert $U_{j}^{m}$ de $\Lambda_{j}$, de mesure inf\'erieure \`a $1/m$, qui contient $\Lambda_{j} \cap \Lambda_{k}$. Comme $U_{j}^{m}$ est la restriction \`a $\Lambda_{j}$ d'un ouvert $V_{j}^{m}$ de $\R^{2n}$, le lemme d'Urysohn nous fournit une fonction $f_{j}^{m} \in C^{\infty}_ {0} (\R^{2n})$, de support inclus dans $V_{j}^{m}$, \'egale \`a $1$ au voisinage de $\Lambda_{j} \cap \Lambda_{k}$. De la m\^eme mani\`ere, on construit $f_{k}^{m}$ en intervertissant les indices $j$ et $k$. On pose alors $f_{m} = 1 - f_{j}^{m} f_{k}^{m}$ que l'on tronque hors du support de $q$.

Soit $\widetilde{q} \in C^{\infty}_{0} (\R^{2n})$ \`a support proche de celui de $q$ et \'egal \`a $1$ pr\`es de celui-ci. D'apr\`es \eqref{a76}, \eqref{a73}, la discussion apr\`es \eqref{a28} et les propri\'et\'es du calcul pseudo-diff\'erentiel, on a
\begin{equation}
\Op (q) u_{h} = \Op (q) (B_{1} + B_{2} ) + R_{0} \quad \text{et} \quad \Op (\widetilde{q} ) u_{h} = \Op ( \widetilde{q} ) (B_{1} + B_{2} ) + \widetilde{R}_{0} ,
\end{equation}
avec $R_{0} , \widetilde{R}_{0} = \CO^{- \beta} (1)$ pour tout $\beta \geq 0$. Le lemme \ref{a7}, \eqref{a24} et \eqref{a54} et la proposition 1.2.2 de \cite{Iv98_01} sur la norme $L^{2}$ d'une int\'egrale oscillante impliquent
\begin{align}
v_{\delta} B_{1} &= R_{1} = \CO ( \sqrt{\delta}h^{- \frac{1}{2}})  \\
\Op (1-f_{m}) (1 - v_{\delta} ) B_{1} &= R_{2} = o_{m\to \infty} (h^{- \frac{1}{2}}) + \CO_{m} (h^{\frac{1}{2}}) \\
\Op (f_{m}) (1 - v_{\delta} ) B_{1} &= h^{-\frac{1}{2}} B_{1}^{m,\delta}  \nonumber \\
&\hspace{-2cm}\text{avec } B_{1}^{m,\delta} \in \dl^{0}_{\text{cl}} ( \widetilde{\Lambda} , (1 -v_{\delta} ) f_{m\vert_{\Lambda_{k}}} b_{0} , \Omega \cap \supp f_{m} (1- v_{\delta}) )  \\
\Op (1-f_{m}) B_{2} &= R_{3} = o_{m\to \infty} (h^{- \frac{1}{2}}) + \CO_{m} (h^{\frac{1}{2}})  \\
\Op (f_{m}) B_{2} &= h^{-\frac{1}{2}} \sum_{k=0}^{K} B_{2}^{k,m} \nonumber  \\
&\hspace{-2cm}\text{avec }  B_{2}^{k,m} \in \dl^{0}_{\text{cl}} ( \Lambda_{k} , f_{m\vert_{\Lambda_{k}}} b^{k}_{0} , \Omega_{k} \cap \supp f_{m})
\end{align}
Ici, une fonction est un $o_{m\to \infty} (1)$ si elle tend vers $0$ lorsque $m \to \infty$. De m\^eme, une fonction est un $\CO_{m} (h)$ si, pour tout $m$ fix\'e, elle est un $\CO (h)$. En particulier,
\begin{equation}
\Op (q) u_{h} = \CO^{- \beta} (h^{- \frac{1}{2}}) ,
\end{equation}
pour tous $\beta \geq 0$ et $q \in C^{\infty}_{0} (\R^{2n})$.

On a alors
\begin{align}
\< \Op (q) u_{h} , u_{h} \> =& \< \Op (q) u_{h} , \Op ( \widetilde{q} ) u_{h} \> + \CO (h^{\infty} )  \nonumber \\
=& \< R_{0} + \Op (q) ( R_{1} + R_{2} + R_{3} ) , \Op ( \widetilde{q} ) u_{h} \> \nonumber  \\
&+ h^{- \frac{1}{2}} \< \Op (q) (B_{1}^{m, \delta} + \sum B_{2}^{k,m} ) , R_{0} + \Op ( \widetilde{q} ) ( R_{1} + R_{2} + R_{3} ) \>  \nonumber  \\
&+ h^{-1} \< \Op (q) (B_{1}^{m, \delta} + \sum B_{2}^{k,m} ) , \Op ( \widetilde{q} ) (B_{1}^{m, \delta} + \sum B_{2}^{k,m} ) \> + \CO (h^{\infty} )   \nonumber  \\
=& h^{-1} \< \Op (q) (B_{1}^{m, \delta} + \sum B_{2}^{k,m} ) , (B_{1}^{m, \delta} + \sum B_{2}^{k,m} ) \>   \nonumber  \\
&+ h^{-1} \big( o_{m \to \infty} (1) + \CO (\sqrt{\delta}) + \CO_{\delta ,m} ( h^{\frac{1}{2}} ) \big) .  \label{a34}
\end{align}
Comme les $K+1$ vari\'et\'es lagrangiennes $\widetilde{\Lambda}$ et $\Lambda_{k}$ ne s'intersectent pas sur $\Omega \cap \supp f_{m} (1- v_{\delta} )$ et $\Omega_{k} \cap \supp f_{m}$, les distributions lagrangiennes $B_{1}^{m, \delta}$ et $B_{2}^{k,m}$ ont des micro-supports deux \`a deux disjoints. Ceci entraine
\begin{gather}
\< \Op (q) B_{1}^{m, \delta} , B_{2}^{k,m} \> = \CO_{m,\delta} (h^{\infty})    \label{a37}   \\
\< \Op (q) B_{2}^{k,m} , B_{2}^{\widetilde{k},m} \> = \CO_{m} (h^{\infty}) \text{ pour } k \neq \widetilde{k} . \label{a38}
\end{gather}

Le calcul de l'application d'un op\'erateur pseudo-diff\'erentiel sur une distribution lagrangienne (th\'eor\`eme 1.2.8 de \cite{Iv98_01}) donne
\begin{equation}  \label{a30}
\Op (q) \dl^{0}_{\text{cl}} ( \Lambda , a_{0} , U )  \subset \dl^{0}_{\text{cl}} ( \Lambda , q_{\vert_{\Lambda}} \times a_{0} , U ) .
\end{equation}
En particulier,
\begin{gather*}
\Op (q) B_{1}^{m,\delta} \in \dl^{0}_{\text{cl}} ( \widetilde{\Lambda} , q (1 -v_{\delta} ) f_{m\vert_{\Lambda_{k}}} b_{0} , \Omega \cap \supp f_{m} (1- v_{\delta}) )  \\
\Op (q) B_{2}^{k,m} \in \dl^{0}_{\text{cl}} ( \Lambda_{k} , q f_{m\vert_{\Lambda_{k}}} b^{k}_{0} , \Omega_{k} \cap \supp f_{m}) .
\end{gather*}
Un calcul \'evident montre alors que
\begin{align}
\< \Op (q) B_{1}^{m, \delta} , B_{1}^{m, \delta} \> =& \int_{\widetilde{\Lambda}}  q \vert (1 - v_{\delta} ) f_{m} b_{0} \vert^{2} \, d \mu + \CO_{m, \delta} (h)   \\
\< \Op (q) B_{2}^{k,m} , B_{2}^{k,m} \> =& \int_{\Lambda_{k}}  q \vert f_{m} b_{0}^{k} \vert^{2} \, d \mu + \CO_{m} (h) ,
\end{align}
o\`u $\mu$ est la mesure de Lebesgue sur $\widetilde{\Lambda}$ ou $\Lambda_{k}$. Comme $b_{0}^{k}$ est r\'eguli\`ere, on a
\begin{equation}
\< \Op (q) B_{2}^{k,m} , B_{2}^{k,m} \> = \int_{\Lambda_{k}} q \vert b_{0}^{k} \vert^{2} \, d \mu + o_{m\to \infty} (1) + \CO_{m} (h) . \label{a35}
\end{equation}
Pour la partie concernant $B_{1}^{m, \delta}$, le symbole $b_{0}$ donn\'e par \eqref{a31} n'est pas r\'egulier en $0$ car $\vert \det \varphi''_{t, \xi} \vert^{\frac{1}{2}}$ explose quand $x \to 0$. Mais en effectuant le changement de variable dont le Jacobien est donn\'e en \eqref{a33} et en utilisant l'estimation \eqref{a32}, on peut \'ecrire
\begin{align}
\< \Op (q) B_{1}^{m, \delta} , B_{1}^{m, \delta} \> =& \int q \vert (1 - v_{\delta} ) f_{m} \vert^{2} \vert b_{0} \vert^{2} \left\vert \frac{d x}{ d (t, \xi )} \right\vert d t \, d \xi + \CO_{m, \delta} (h)  \nonumber  \\
=& \int q \vert b_{0} \vert^{2} (x (t, \xi), \nabla \psi (x (t, \xi ))) \left\vert \frac{d x}{ d (t, \xi )} \right\vert d t \, d \xi   \nonumber  \\
&+ o_{m \to \infty}(1) + o_{\delta \to 0} (1) + \CO_{m, \delta} (h)  . \label{a36}
\end{align}
Dans la formule pr\'ec\'edente, l'int\'egration se fait sur $[0, \varepsilon_{0} [ \times \sqrt{2(E - V(0))}\S^{n-1}$ pour $\varepsilon_{0}$ assez petit. On reconnait une mesure dans le premier terme de \eqref{a36}. La proposition r\'esulte alors de \eqref{a34}, \eqref{a37}, \eqref{a38}, \eqref{a35} et \eqref{a36} : on fixe d'abord $\delta$ assez petit et $m$ assez grand et finalement, on fait tendre $h$ vers $0$.

Ceci conclut la preuve pour $E = E_{0} + h F_{1}$. Dans le cas $E = E_{0} + h E_{1} + o_{h \to 0} (h)$, il r\'esulte des constructions pr\'ec\'edentes que les symboles $b_{0}$ et $b_{0}^{k}$ qui apparaissent dans \eqref{a35} et \eqref{a36} sont holomorphes par rapport a $F_{1}$ et que les termes d'erreur sont localement uniformes par rapport \`a $F_{1}$. Il suffit alors de prendre $F_{1} = E_{1} + o_{h \to 0} (1)$ et d'utiliser la formule de Taylor en $F_{1} = E_{1}$.
\end{proof}

\Subsection{Calcul de la mesure $\mu_{h^{-1}}$.}

Dans cette partie on calcule la mesure de Wigner $\mu_{h^{-1}}$ d\'efinit dans la proposition \ref{a39}. On a la proposition suivante

\begin{prop}\sl \label{a41}
La mesure $\mu_{h^{-1}}$ poss\`ede les propri\'et\'es suivantes qui la caract\'erisent uniquement :
\begin{enumerate}
\item La mesure $\mu_{h^{-1}}$ est support\'ee dans la surface d'\'energie $E_{0}$ :
\begin{equation} \label{a44}
\supp \mu_{h^{-1}} \subset p^{-1}(E_{0})
\end{equation}

\item La mesure $\mu_{h^{-1}}$ est nulle dans la zone entrante :
\begin{equation} \label{a40}
\mu_{h^{-1}} =0 \text{ dans } \{ (x, \xi ) ; \ \vert x \vert \text{ grand  et } \cos (x, \xi ) < -1/2 \}.
\end{equation}

\item La mesure $\mu_{h^{-1}}$ v\'erifie l'\'equation d'\'evolution :
\begin{align}
( {\rm H}_{p} + 2 \im E_{1} ) \, \mu_{h^{-1}} =&  \big( \xi \partial_{x} - \nabla V (x) \partial_{\xi} + 2 \im E_{1} \big) \mu_{h^{-1}}  \nonumber  \\
=& ( 2 \pi )^{1-n} \vert \widehat{S} (\xi ) \vert^{2} \delta_{x =0} \delta_{\frac{1}{2} \xi^{2} + V(0) = E_{0}}  \label{a42}
\end{align}
\end{enumerate}
\end{prop}

\begin{proof}
L'estimation \eqref{a45} combin\'ee avec l'existence de la mesure de Wigner (proposition \ref{a39}) implique \eqref{a44}.

La preuve de \eqref{a40} r\'esulte de l'existence de la mesure et du fait que $u_{h} = \CO^{\alpha} (h^{\infty})$ dans la zone entrante (\'equation \eqref{a77} et proposition \ref{a43}).

Soit $q \in C^{\infty}_{0} (\R^{2n})$ localis\'e pr\`es de $\rho_{0} = (x_{0} , \xi_{0} )$. Comme $(P - E) u_{h} = S_{h}$, on a $(P - E) u_{h} = \CO (h^{\infty})$ micro-localement pr\`es de $\rho_{0}$, si $x_{0} \neq 0$. D'apr\`es \ref{h6}, $E - \overline{E} = 2 i \im E_{1} h + o (h)$. En remarquant que ${\rm H}_{p} (q)$ est le symbole principal de l'op\'erateur pseudo-diff\'erentiel $i h^{-1} [ P ,  \Op (q) ]$, il vient
\begin{align}
\< q , ( {\rm H}_{p} + 2 \im E_{1} ) \, \mu_{h^{-1}} \> =& - \< {\rm H}_{p} (q) - 2 q \im E_{1} , \mu_{h^{-1}} \>  \nonumber  \\
=& - \lim_{h \to 0} \< i \big( (P - E)^{*} \Op (q) - \Op (q) (P - E) \big) u_{h} ,u_{h} \> + o(1)=0 .    \label{a46}
\end{align}
Ainsi, $( {\rm H}_{p} + 2 \im E_{1} ) \, \mu_{h^{-1}} =0$ pour $x \neq 0$.

On suppose maintenant que $x_{0} =0$. Pr\`es de $\rho_{0}$, la mesure $\mu_{h^{-1}}$ est donn\'ee comme somme des termes dominants de \eqref{a35} et \eqref{a36}. D'apr\`es \eqref{a71}, $B_{2}^{k}$ v\'erifie $(P - E) B_{2}^{k} =0$ micro-localement pr\`es de $\rho_{0}$ et la preuve de \eqref{a46} montre que ces termes ne contribuent pas dans le calcul de $( {\rm H}_{p} + 2 \im E_{1} ) \, \mu_{h^{-1}}$. Donc,
\begin{align}
\< q , ( {\rm H}_{p} + 2 \im E_{1} ) \, \mu_{h^{-1}} \> =& - \int \big( {\rm H}_{p} (q) - 2 q \im E_{1} \big) \vert b_{0} \vert^{2} (x (t, \xi), \nabla \psi (x (t, \xi ))) \left\vert \frac{d x}{ d (t, \xi )} \right\vert d t \, d \xi  \nonumber  \\
=& - \int \big( \partial_{t}  q (\rho (t , \xi )) - 2 q (\rho (t , \xi )) \im E_{1} \big) \vert b_{0} \vert^{2} (x (t, \xi)) \left\vert \frac{d x}{ d (t, \xi )} \right\vert d t \, d \xi ,
\end{align}
o\`u $\rho (t, \xi )$ est la courbe hamiltonienne avec donn\'ee initiale $( 0 , \xi )$ en $t=0$. Comme avant, l'int\'egration se fait sur $[0, \varepsilon_{0} [ \times \sqrt{2(E - V(0))}\S^{n-1}$. D'apr\`es \eqref{a47}, \eqref{a31}, \eqref{a32} et \eqref{a33}, on~a
\begin{equation}
\vert b_{0} \vert^{2} (x (t, \xi)) \left\vert \frac{d x}{ d (t, \xi )} \right\vert = (2 \pi )^{1-n} \vert \widehat{S} (\xi) \vert^{2} (1 + o_{t \to 0} (1)) .
\end{equation}
Soit $\chi \in C^{\infty}_{0} (\R , [0,1])$ avec $\chi =1$ pr\`es de $0$ et on note $\chi_{m} (t) = \chi (t m)$. On voit $\chi_{m}$ comme une fonction d\'efinie sur un voisinage de $\rho_{0} \in \R^{2n}$ en posant $\chi_{m} ( x , \eta) = \chi_{m} (t)$ si $x = x (t , \xi)$ (rappelons que $(t , \xi ) \to x ( t, \xi )$ est le changement de variables d\'efini avant \eqref{a33}). En utilisant la discussion pr\'ec\'edente, il vient
\begin{align}
\< q , ( {\rm H}_{p} + 2 \im & E_{1} ) \, \mu_{h^{-1}} \>   \nonumber  \\
=& \< \chi_{m} q , ( {\rm H}_{p} + 2 \im E_{1} ) \, \mu_{h^{-1}} \> + \< (1 - \chi_{m} ) q , ( {\rm H}_{p} + 2 \im E_{1} ) \, \mu_{h^{-1}} \>  \nonumber  \\
=& - \< ( {\rm H}_{p} - 2 \im E_{1} ) \, ( \chi_{m} q) , \mu_{h^{-1}} \> + 0   \nonumber  \\
=& - \iint_{0}^{\infty} \big( \chi (t m) ( \partial_{t} q - 2 q \im E_{1} )+ m \chi' (t m) q \big) (2 \pi )^{1-n} \vert \widehat{S} (\xi) \vert^{2} (1 + o_{t \to 0} (1)) \, d t \, d \xi \nonumber  \\
=& \CO (m^{-1} ) - \iint_{0}^{\infty} m \chi' (t m) q ( \rho (0, \xi ) ) (2 \pi )^{1-n} \vert \widehat{S} (\xi) \vert^{2} (1 + o_{t \to 0} (1)) \, d t \, d \xi \nonumber   \\
=& \int  q (0 , \xi ) (2 \pi )^{1-n} \vert \widehat{S} (\xi) \vert^{2} \, d \xi + o_{m \to \infty} (1)  \nonumber  \\
=& \int  q (0 , \xi ) (2 \pi )^{1-n} \vert \widehat{S} (\xi) \vert^{2} \, d \xi ,
\end{align}
o\`u l'int\'egration se fait sur $\sqrt{2(E - V(0))}\S^{n-1}$.
\end{proof}

\section{Contre-exemple quand l'hypoth\`ese \ref{h5} n'est pas v\'erifi\'ee.}
\label{a50}

On donne ici un exemple (en dimension $n=1$ puis en dimension quelconque) de potentiel et de fonction $S$ qui ne v\'erifient pas \ref{h5} et dont la mesure limite n'est pas unique. Pour fixer les id\'ees, on suppose que $E =E_{0}$.

En dimension $n=1$, consid\'erons un potentiel satisfaisant les hypoth\`eses \ref{h1}, \ref{h2} et \ref{h3} de la forme suivante :
\begin{figure}[!h]
\begin{center}
\begin{picture}(0,0)%
\includegraphics{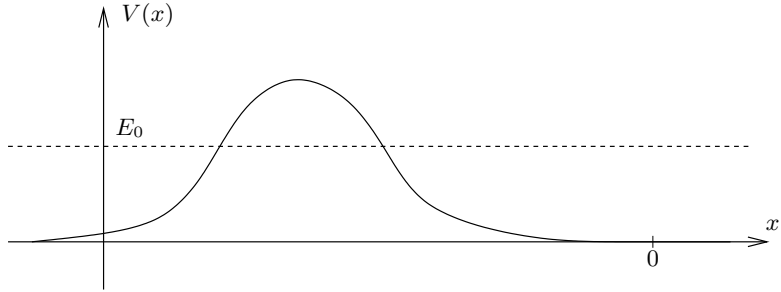}%
\end{picture}%
\setlength{\unitlength}{1973sp}%
\begingroup\makeatletter\ifx\SetFigFont\undefined%
\gdef\SetFigFont#1#2#3#4#5{%
  \reset@font\fontsize{#1}{#2pt}%
  \fontfamily{#3}\fontseries{#4}\fontshape{#5}%
  \selectfont}%
\fi\endgroup%
\begin{picture}(9644,3644)(1179,-3383)
\put(2551,-1411){\makebox(0,0)[lb]{\smash{{\SetFigFont{9}{10.8}{\rmdefault}{\mddefault}{\updefault}$E_{0}$}}}}
\put(2626, 14){\makebox(0,0)[lb]{\smash{{\SetFigFont{9}{10.8}{\rmdefault}{\mddefault}{\updefault}$V(x)$}}}}
\put(10726,-2611){\makebox(0,0)[lb]{\smash{{\SetFigFont{9}{10.8}{\rmdefault}{\mddefault}{\updefault}$x$}}}}
\put(9226,-3061){\makebox(0,0)[lb]{\smash{{\SetFigFont{9}{10.8}{\rmdefault}{\mddefault}{\updefault}$0$}}}}
\end{picture}%
\caption{Le potentiel $V(x)$.} \label{a57}
\end{center}
\end{figure}\newline
L'ensemble $\Lambda$ ressemble alors \`a :
\begin{figure}[!h]
\begin{center}
\begin{picture}(0,0)%
\includegraphics{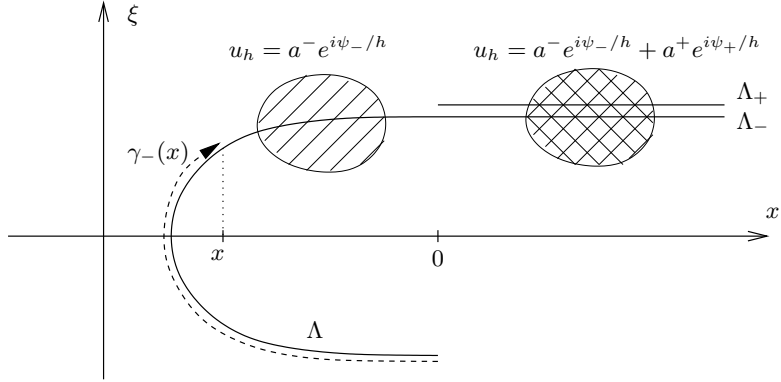}%
\end{picture}%
\setlength{\unitlength}{1973sp}%
\begingroup\makeatletter\ifx\SetFigFont\undefined%
\gdef\SetFigFont#1#2#3#4#5{%
  \reset@font\fontsize{#1}{#2pt}%
  \fontfamily{#3}\fontseries{#4}\fontshape{#5}%
  \selectfont}%
\fi\endgroup%
\begin{picture}(9644,4844)(1179,-3383)
\put(6526,-1936){\makebox(0,0)[lb]{\smash{{\SetFigFont{9}{10.8}{\rmdefault}{\mddefault}{\updefault}$0$}}}}
\put(10726,-1336){\makebox(0,0)[lb]{\smash{{\SetFigFont{9}{10.8}{\rmdefault}{\mddefault}{\updefault}$x$}}}}
\put(2701,1139){\makebox(0,0)[lb]{\smash{{\SetFigFont{9}{10.8}{\rmdefault}{\mddefault}{\updefault}$\xi$}}}}
\put(3751,-1861){\makebox(0,0)[lb]{\smash{{\SetFigFont{9}{10.8}{\rmdefault}{\mddefault}{\updefault}$x$}}}}
\put(2701,-586){\makebox(0,0)[lb]{\smash{{\SetFigFont{9}{10.8}{\rmdefault}{\mddefault}{\updefault}$\gamma_{-} (x)$}}}}
\put(3976,689){\makebox(0,0)[lb]{\smash{{\SetFigFont{9}{10.8}{\rmdefault}{\mddefault}{\updefault}$u_{h} =a^{-} e^{i \psi_{-} /h}$}}}}
\put(7051,689){\makebox(0,0)[lb]{\smash{{\SetFigFont{9}{10.8}{\rmdefault}{\mddefault}{\updefault}$u_{h} =a^{-} e^{i \psi_{-} /h} + a^{+} e^{i \psi_{+} /h}$}}}}
\put(4951,-2836){\makebox(0,0)[lb]{\smash{{\SetFigFont{9}{10.8}{\rmdefault}{\mddefault}{\updefault}$\Lambda$}}}}
\put(10351,-211){\makebox(0,0)[lb]{\smash{{\SetFigFont{9}{10.8}{\rmdefault}{\mddefault}{\updefault}$\Lambda_{-}$}}}}
\put(10351,164){\makebox(0,0)[lb]{\smash{{\SetFigFont{9}{10.8}{\rmdefault}{\mddefault}{\updefault}$\Lambda_{+}$}}}}
\end{picture}%
\caption{La vari\'et\'e $\Lambda$ avec $\Lambda_{-} = \Lambda_{+}$ dans la zone quadrill\'ee.} \label{a58}
\end{center}
\end{figure}

Soit $S$ une fonction v\'erifiant \ref{h4}. D'apr\`es la partie \ref{a55}, la solution $u_{h}$, micro-localement dans la zone hachur\'ee, est de la forme
\begin{equation}
u_{h} = a^{-}( x,h) e^{i \psi_{-} (x) /h} + \CO (h^{\infty}) ,
\end{equation}
o\`u $a^{-}$ est un symbole classique et $\psi = \int_{\gamma_{-} (x)} \xi \, d x = \int_{0}^{t (x)} \xi (s)^{2} \, d s$. En particulier, si le support de $q \in C^{\infty}_{0} (\R^{2n})$ est contenu dans la zone hachur\'ee,
\begin{equation}
h \< \Op (q) u_{h} , u_{h} \> \longrightarrow \< q , \mu_{h^{-1}} \> , \text{ pour } h \to 0 ,
\end{equation}
et $\mu_{h^{-1}}$ v\'erifie les conclusions du th\'eor\`eme \ref{a48}.

Par contre, micro-localement dans la zone quadrill\'ee, $u_{h}$ est de la forme
\begin{equation}
u_{h} = a^{-}( x,h) e^{i \psi_{-} (x) /h} + a^{+}( x,h) e^{i \psi_{+} (x) /h} + \CO (h^{\infty}) ,
\end{equation}
o\`u $a^{+}$ et $\psi_{+}$ v\'erifient le m\^eme type de propri\'et\'es que $a^{-}$ et $\psi_{-}$. En particulier,
\begin{equation}
\psi_{-} (x) = \int_{\gamma_{-} (x)} \xi \, d x = \int_{\gamma_{-} (0)} \xi \, d x + \int_{\gamma_{+} (x)} \xi \, d x = {\mathcal A} + \psi_{+} (x) ,
\end{equation}
o\`u ${\mathcal A} = \int_{0}^{t (0)} \xi (s)^{2} \, d s >0$ est l'action associ\'ee \`a la courbe $\gamma_{-} (0)$. En outre, en utilisant \eqref{a31}, le fait que les $a^{\bullet}$ v\'erifient la m\^eme \'equation de transport dans (la projection en $x$ de) la zone quadrill\'ee et la th\'eorie de Maslov, on trouve
\begin{equation} \label{a65}
a^{-}_{0} (x) = e^{- i \nu \pi /2} \frac{\widehat{S} \big( - \sqrt{2 (E-V(0))} \big)}{\widehat{S} \big( \sqrt{2 (E-V(0))} \big)} r a^{+}_{0} (x) = e^{i \theta} \widetilde{r} a^{+}_{0} (x) ,
\end{equation}
o\`u $\nu$ est l'indice de Maslov associ\'e \`a la courbe $\gamma_{-} (0)$ et $r, \widetilde{r} \in ] 0 , + \infty [$, $\theta \in \R$ sont des constantes. Pour que $a^{-}_{0}$ et $a^{+}_{0}$ soient non nuls, on suppose que les transform\'ees de Fourier de $S$ qui apparaissent dans \eqref{a65} sont non nulles. Ainsi, si le support de $q \in C^{\infty}_{0} (\R^{2n})$ est contenu dans la zone quadrill\'ee,
\begin{align}
h \< \Op (q) u_{h} , u_{h} \> =& \int q (x, \nabla \psi_{\bullet} (x)) \big( \vert a^{-}_{0} \vert^{2} + \vert a^{+}_{0} \vert^{2} + 2 \re \big( a^{-}_{0} \overline{a^{+}_{0}} e^{i {\mathcal A} /h} \big) \big) d x + o_{h \to 0} (1)   \nonumber  \\
=& \int q (x, \nabla \psi_{\bullet} (x)) \big( \vert a^{-}_{0} \vert^{2} + \vert a^{+}_{0} \vert^{2} + 2 \vert a^{-}_{0} \vert \vert a^{+}_{0} \vert \cos (\theta + {\mathcal A} /h) \big) d x + o_{h \to 0} (1) .
\end{align}
Comme ${\mathcal A} \neq 0$, on obtient facilement la

\begin{prop}\sl  \label{a51}
Sous les hypoth\`eses pr\'ec\'edentes et pour tout $\nu \in [-1 , 1]$, on peut trouver un ensemble de $h$ adh\'erant \`a $0$ tel que
\begin{equation}
h \< \Op (q) u_{h} , u_{h} \> \longrightarrow \int q (x, \nabla \psi_{\bullet} (x)) \big( \vert a^{-}_{0} \vert^{2} + \vert a^{+}_{0} \vert^{2} + 2 \nu \vert a^{-}_{0} \vert \vert a^{+}_{0} \vert \big) \, d x ,
\end{equation}
quand $h$ tend vers $0$ dans cet ensemble.
\end{prop}

\begin{figure}
\begin{center}
\hspace{-2cm}\includegraphics[scale=0.4]{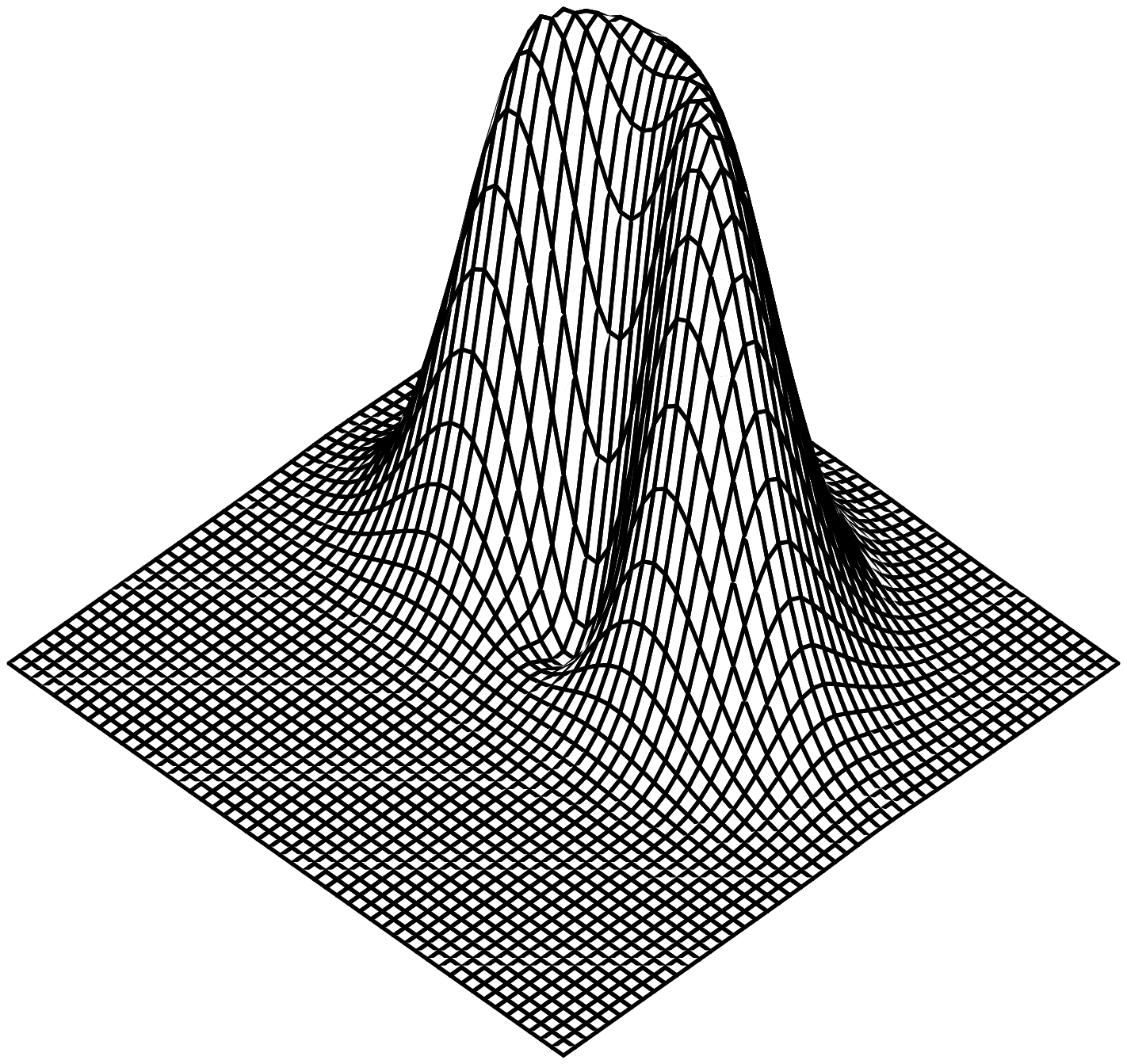}%
\begin{picture}(0,0)%
\includegraphics{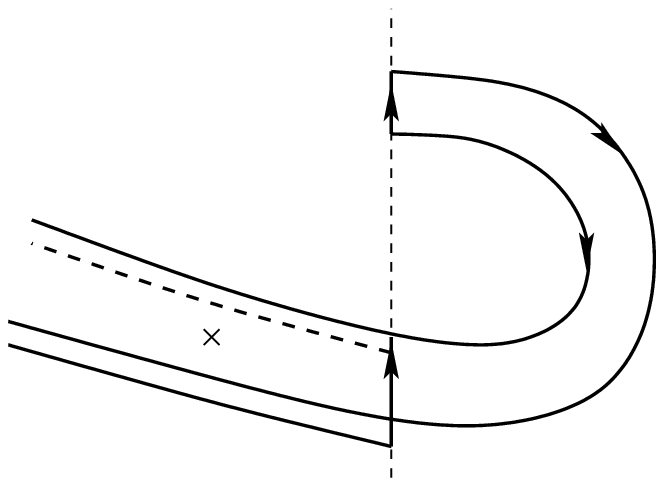}%
\end{picture}%
\setlength{\unitlength}{1973sp}%
\begingroup\makeatletter\ifx\SetFigFont\undefined%
\gdef\SetFigFont#1#2#3#4#5{%
  \reset@font\fontsize{#1}{#2pt}%
  \fontfamily{#3}\fontseries{#4}\fontshape{#5}%
  \selectfont}%
\fi\endgroup%
\begin{picture}(6705,4984)(1861,-4723)
\put(8176,-886){\makebox(0,0)[lb]{\smash{{\SetFigFont{9}{10.8}{\rmdefault}{\mddefault}{\updefault}$\gamma_{-,x_{1}} (0)$}}}}
\put(5776,-4636){\makebox(0,0)[lb]{\smash{{\SetFigFont{9}{10.8}{\rmdefault}{\mddefault}{\updefault}$x=0$}}}}
\put(5626,-736){\makebox(0,0)[lb]{\smash{{\SetFigFont{9}{10.8}{\rmdefault}{\mddefault}{\updefault}$l_{1}$}}}}
\put(6151,-3361){\makebox(0,0)[lb]{\smash{{\SetFigFont{9}{10.8}{\rmdefault}{\mddefault}{\updefault}$l_{2}$}}}}
\put(4501,-3061){\makebox(0,0)[lb]{\smash{{\SetFigFont{9}{10.8}{\rmdefault}{\mddefault}{\updefault}$\rho_{0}$}}}}
\put(8551,-3361){\makebox(0,0)[lb]{\smash{{\SetFigFont{9}{10.8}{\rmdefault}{\mddefault}{\updefault}$\Lambda$}}}}
\put(1876,-3136){\makebox(0,0)[lb]{\smash{{\SetFigFont{9}{10.8}{\rmdefault}{\mddefault}{\updefault}$\Lambda_{+}$}}}}
\put(1876,-2386){\makebox(0,0)[lb]{\smash{{\SetFigFont{9}{10.8}{\rmdefault}{\mddefault}{\updefault}$\Lambda_{-}$}}}}
\put(7201,-3361){\makebox(0,0)[lb]{\smash{{\SetFigFont{9}{10.8}{\rmdefault}{\mddefault}{\updefault}$\Sigma$}}}}
\put(6676,-2161){\makebox(0,0)[lb]{\smash{{\SetFigFont{9}{10.8}{\rmdefault}{\mddefault}{\updefault}$\gamma_{-,x_{2}} (0)$}}}}
\end{picture}%
\caption{Un exemple de potentiel, \'egal \`a l'oscillateur harmonique dans un petit secteur angulaire, qui re-focalise en $0$ tous les rayons issus d'un secteur angulaire en $0$. L'ensemble $\Lambda$, avec $\Lambda_{-} = \Lambda_{+}$ pr\`es de $\rho_{0}$, est d\'ecrit \`a droite.} \label{a56}
\end{center}
\end{figure}

En particulier, la mesure limite n'est plus unique. La mesure du th\'eor\`eme \ref{a48} correspond \`a~$\nu =0$. Il se peut que la mesure soient nulle en dehors d'un compact. En dimension plus grande, en s'inspirant de \cite{Ca05_01}, on peut construire un potentiel tel que localement pr\`es d'un point $\rho_{0}$ l'ensemble $\Lambda$ soit la r\'eunion de 2 vari\'et\'es $\Lambda_{-}$ et $\Lambda_{+}$ qui coïncident (figure \ref{a56}).

Comme en dimension $n=1$, micro-localement pr\`es de $\rho_{0}$, la fonction $u_{h}$ s'\'ecrit
\begin{equation}
u_{h} = a^{-}( x,h) e^{i \psi_{-} (x) /h} + a^{+}( x,h) e^{i \psi_{+} (x) /h} + \CO (h^{\infty}) ,
\end{equation}
La diff\'erence des phases ${\mathcal A} (x) = \psi_{-} (x) - \psi_{+} (x)$ est encore l'action le long de la courbe $\gamma_{- , x} (0)$ :
\begin{equation}
{\mathcal A} (x) = \int_{\gamma_{- ,x} (0)} \xi \cdot d x .
\end{equation}
Mais on remarque que cette action ne d\'epend pas de $x$. En effet, si $\gamma_{- , x_{1}} (0)$ et $\gamma_{- , x_{2}} (0)$ sont deux chemins, le fait que les $l_{\bullet}$ soient dans $\{ x = 0 \}$, la formule de Stokes et le caract\`ere lagrangien de $\Lambda$ impliquent
\begin{align}
\int_{\gamma_{- , x_{1}} (0)} \xi \cdot d x - \int_{\gamma_{- , x_{2}} (0)} \xi \cdot d x =& \int_{\gamma_{- , x_{1}} (0)} \xi \cdot d x + \int_{l_{2}} \xi \cdot d x - \int_{\gamma_{- , x_{2}} (0)} \xi \cdot d x + \int_{l_{1}} \xi \cdot d x   \nonumber  \\
=& \int_{\Sigma} d \xi \wedge d x  = 0.  \label{a52}
\end{align}
On peut aussi dire que comme $\Lambda_{-}$ et $\Lambda_{+}$ coïncident, on a $\nabla \psi_{-} = \nabla \psi_{+}$ et donc $\psi_{-}$ et $\psi_{+}$ diff\`erent d'une constante. Cette remarque permet d'adapter la preuve de la proposition \ref{a51} \`a cette situation.

\section{Preuve de la proposition \ref{a59}.}
\label{a60}

En notant $a \# b$ le symbole de l'op\'erateur $\Op a \circ \Op b$, on a
\begin{align}
q (p-E_{0})^{2} =& (p-E_{0}) \big( q \# (p-E_{0}) \big) + \frac{h}{2 i} (p-E_{0}) {\rm H}_{p} \, q + h^{2} r (x, \xi ,h)  \nonumber  \\
=& (p-E_{0}) \# \big( q \# (p-E_{0}) \big) - \frac{h}{2 i} {\rm H}_{p} \, ( q \# (p-E_{0}) ) + \frac{h}{2 i} (p-E_{0}) {\rm H}_{p} \, q + h^{2} r (x, \xi ,h)   \nonumber  \\
=& (p-E_{0}) \# \big( q \# (p-E_{0}) \big) + h^{2} r (x, \xi ,h) ,
\end{align}
o\`u $r \in S_{{\rm cl}} (1)$ change de lignes en lignes et \`a son support inclus dans le support (compact) de $q$, modulo $S (h^{\infty} \< x \>^{- \infty} \< \xi \>^{- \infty} )$. Mais alors, en utilisant que $u_{h} = \CO (h^{-1/2})$ micro-localement pr\`es du support de $q$, $\Vert S_{h} \Vert = \CO (1)$ et $\< \Op ( q ) u_{h} , u_{h} \> = \CO (h^{-1})$ (d'apr\`es le th\'eor\`eme \ref{a48}), il vient
\begin{align}
\< \Op ( q (p-E_{0} )^{2} ) u_{h} , u_{h} \> =& \< \Op ( p-E_{0} ) \Op ( q ) \Op ( p-E_{0} ) u_{h} , u_{h} \> + \CO (h)  \nonumber  \\
=& \< \Op ( q ) \Op ( p-E ) u_{h} , \Op ( p-E ) u_{h} \> + \CO (h) \Vert \Op ( q ) u_{h} \Vert \Vert S_{h} \Vert + \CO (h)   \nonumber   \\
=& \< \Op ( q ) S_{h} , S_{h} \> + \CO (h^{\frac{1}{2}}) .  \label{a61}
\end{align}
Grâce \`a la formule sur l'application d'un op\'erateur pseudo-diff\'erentiel sur une distribution lagrangienne (th\'eor\`eme 1.2.8 de \cite{Iv98_01}) ou en faisant le calcul explicitement, on obtient
\begin{align}
\Op ( q ) S_{h} =& \Op ( q ) g(x) S_{h} + \CO (h^{\infty}) \nonumber  \\
=& (2 \pi )^{-n} h^{-n/2} \int \Op ( q ) g (x) e^{i x \cdot \xi / h} \widehat{S} ( \xi) \, d \xi + \CO (h^{\infty})   \nonumber  \\
=& (2 \pi )^{-n} h^{-n/2} \int q (x, \xi ) g (x) e^{i x \cdot \xi / h} \widehat{S} ( \xi) \, d \xi + \CO (h) ,
\end{align}
pour $g (x) \in C^{\infty}_{0} (\R^{n})$ avec $g =1$ pr\`es de $0$. En particulier, \eqref{a61} devient
\begin{equation*}
\< \Op ( q (p-E_{0} )^{2} ) u_{h} , u_{h} \> = (2 \pi )^{-2 n} h^{-n} \iiint q (x, \xi ) g (x) e^{i x \cdot \xi / h} \widehat{S} ( \xi) e^{- i x \cdot \theta /h} \overline{\widehat{S}} ( \theta ) \, d x \, d \xi \, d \theta + \CO (h) .
\end{equation*}
La m\'ethode de la phase stationnaire en variables $( x, \theta )$ donne alors
\begin{equation}
\< \Op ( q (p-E_{0} )^{2} ) u_{h} , u_{h} \> = (2 \pi )^{- n} \int q (0 , \xi ) \vert \widehat{S} ( \xi) \vert^{2} \, d \xi + \CO (h) ,
\end{equation}
ce qui d\'emontre la proposition.

\section{Equivalence entre les hypoth\`eses \ref{h5} et \ref{h7}.}
\label{a85}

\begin{lemm}\sl \label{a84}
Soient $\Lambda_{1}$ et $\Lambda_{2}$ deux sous-vari\'et\'es de $\R^{n}$ de m\^eme dimension $m$. Alors l'ensemble ferm\'e  $\Lambda_{1} \cap \Lambda_{2}$ est de mesure nulle dans $\Lambda_{1}$ si et seulement si il est de mesure nulle dans $\Lambda_{2}$.
\end{lemm}

\begin{proof}
On travaille localement pr\`es de $\rho_{0} \in \Lambda_{1} \cap \Lambda_{2}$. Il est possible de trouver, $E$, un sous-espace vectoriel de $\R^{n}$ de dimension $m-n$ qui soit suppl\'ementaire \`a ${\rm T}_{\rho_{0}} \Lambda_{1}$ et \`a ${\rm T}_{\rho_{0}} \Lambda_{2}$. Pour cela, on peut utiliser la propri\'et\'e suivante : si $A , B$ sont deux sous-espaces vectoriels d'un espace vectoriel $G$, alors $A \cup B = G$ si est seulement si $A=G$ ou $B=G$.

On note $F= {\rm T}_{\rho_{0}} \Lambda_{1}$ et $\Pi$ la projection lin\'eaire de $\R^{n}$ sur $F$ parall\`element \`a $E$. Alors, $\Pi$ est un diff\'eomorphisme local de $\Lambda_{1}$ sur $F$ et de $\Lambda_{2}$ sur $F$. De plus, les mesures de Lebesgue sur $\Lambda_{\bullet}$ sont transport\'ees en des mesures sur $F$. Ces deux mesures ainsi que la mesure de Lebesgue sur $F$ sont absolument continues les unes par rapport aux autres. En particulier, $\Lambda_{1} \cap \Lambda_{2}$ est de mesure nulle dans $\Lambda_{1}$ ssi $\Pi (\Lambda_{1} \cap \Lambda_{2})$ est de mesure nulle dans $F$ pour une des trois mesures pr\'ec\'edentes ssi $\Lambda_{1} \cap \Lambda_{2}$ est de mesure nulle dans $\Lambda_{2}$.
\end{proof}

On rappelle que si $\rho (t) = ( x(t) , \xi (t))$ est une courbe hamiltonienne dans $p^{-1} (E_{0})$, alors $\widetilde{\rho} (t) = ( x(-t) , - \xi (-t))$ en est aussi une. Donc,
\begin{align*}
\big\{ \xi \in \sqrt{2(E_{0} - V(0))} & \S^{n-1} ; \ \exists t < 0 \quad \Pi_{x} \exp ( t {\rm H}_{p}) (0 , \xi) = 0 \big\}    \\
= - & \big\{ \xi \in \sqrt{2(E_{0} - V(0))} \S^{n-1} ; \ \exists t > 0 \quad \Pi_{x} \exp ( t {\rm H}_{p}) (0 , \xi) = 0 \big\} .
\end{align*}
En particulier, dans l'hypoth\`ese \ref{h5}, on peut remplacer $t>0$ par $t \neq 0$.

\begin{lemm}
Supposons \ref{h1}--\ref{h3}. Les hypoth\`eses \ref{h5} et \ref{h7} sont \'equivalentes.
\end{lemm}

\noindent
L'\'equivalence entre les hypoth\`eses \ref{h8} et \ref{h9} se d\'emontre de la m\^eme fa\c{c}on.

\begin{proof}
Supposons \ref{h5} et consid\'erons l'intersection de $\Lambda_{j}$ et $\Lambda_{k}$ avec $j < k$ pr\`es d'un point $\rho_{0}$. La vari\'et\'e $\Lambda_{j}$ est param\'etr\'ee par l'application $\rho (t, \xi) = \exp (t {\rm H}_{p} ) ( 0 , \xi)$ avec $t$ au voisinage de $t_{j}$ et $\xi \in \sqrt{2(E_{0} - V(0))} \S^{n-1}$ au voisinage de $\xi_{j}$.

Donc, un point $\rho = \rho (t, \xi)$ est dans $\Lambda_{j} \cap \Lambda_{k}$ ssi il existe $s$ proche de $t_{k}$ tel que $\rho (s, \xi) = \rho (t, \xi )$ ssi il existe $u = s-t$ proche de $t_{k} - t_{j}$ tel que $\rho (u, \xi ) = (0 , \xi )$ ssi
\begin{align*}
\xi \in M_{j ,k} = \big\{ \xi \in \sqrt{2(E_{0} - V(0))} & \S^{n-1} \text{ proche de } \xi_{j} ;   \\
\exists t > 0 & \text{ proche de } t_{k} - t_{j} \quad \Pi_{x} \exp ( t {\rm H}_{p}) (0 , \xi) = 0 \big\} .
\end{align*}
En utilisant le th\'eor\`eme de Fubini pour s\'eparer les variables $t$ et $\xi$, on a
\begin{align}
{\rm mes}_{n} \, \Lambda_{j} \cap \Lambda_{k} =& \int_{\Lambda_{j}} {\bf 1}_{\Lambda_{j} \cap \Lambda_{k}} d \mu  \lesssim \iint {\bf 1}_{\rho (t, \xi ) \in \Lambda_{j}} {\bf 1}_{\xi \in M_{j ,k}} d t \, d \xi \nonumber  \\
\lesssim& \int {\bf 1}_{\xi \in M_{j ,k}} d \xi  = {\rm mes}_{n-1} \, M_{j ,k} .  \label{a86}
\end{align}
Comme ${\rm mes}_{n-1} M_{j ,k} =0$ d'apr\`es \ref{h5}, l'hypoth\`ese \ref{h7} est v\'erifi\'ee.

R\'eciproquement, supposons \ref{h7}. Soit $\xi_{0} \in \sqrt{2(E_{0} - V(0))} \S^{n-1}$. En utilisant le m\^eme raisonnement que \eqref{a86} mais en sens inverse, on trouve
\begin{align*}
&{\rm mes}_{n-1} \big\{ \xi \in \sqrt{2(E_{0} - V(0))} \S^{n-1} \text{ proche de } \xi_{0}; \ \exists t > 0 \quad \Pi_{x} \exp ( t {\rm H}_{p}) (0 , \xi) = 0 \big\}  \nonumber  \\
&= \sum_{k} M_{0 ,k} \lesssim \sum_{k} {\rm mes}_{n} \, \widetilde{\Lambda} \cap \Lambda_{k} .
\end{align*}
D'apr\`es l'hypoth\`ese \ref{h7}, cette quantit\'e est nulle. L'hypoth\`ese \ref{h5} d\'ecoule alors d'un argument de compacit\'e en $\xi_{0}$.
\end{proof}

\bibliographystyle{amsplain}

\providecommand{\MRhref}[2]{%
  \href{http://www.ams.org/mathscinet-getitem?mr=#1}{#2}
}
\providecommand{\href}[2]{#2}


\end{document}